%% file: swepmor.tex
\documentclass[10pt,a4paper]{article}

\usepackage{graphicx}        % standard LaTeX graphics tool
\graphicspath{{./fig/}}
\usepackage{times}
\usepackage{amsmath,amssymb}
\usepackage{amsthm}
\usepackage{mathrsfs}
\usepackage[a4paper,margin=0.7in,footskip=0.25in]{geometry}

\usepackage{latexsym,bm,multirow}        % for things like \Box and \Diamond
\usepackage{float}
\usepackage{placeins}

\usepackage{caption}
\usepackage{subcaption}
\usepackage{hyperref}
\usepackage{cleveref}
\usepackage{xcolor}
\usepackage[software]{mymacros}
\usepackage{bbm}
\usepackage{tikz}
\usepackage{pgfplots}
\usepgfplotslibrary{colormaps}%
\usepackage{nicefrac}
\usetikzlibrary{external}

\def\addlegendimage{\csname pgfplots@addlegendimage\endcsname}

\usetikzlibrary{external}
\usetikzlibrary{calc,patterns,
                 decorations.pathmorphing,
                 decorations.markings}
\newlength\fheight
\newlength\fwidth

\begin{document}

\title{Data-Driven Learning of Reduced-order Dynamics for a Parametrized Shallow Water Equation}

\author{S\"uleyman Y{\i}ld{\i}z\thanks{Institute of Applied Mathematics, Middle East Technical University, 06800 Ankara, Turkey {yildiz.suleyman@metu.edu.tr}}, Pawan Goyal\thanks{Max Planck Institute for Dynamics of Complex Technical Systems, Sandtorstra{\ss}e 1, D-39106 Magdeburg, Germany {goyalp@mpi-magdeburg.mpg.de}}, Peter Benner\thanks{Max Planck Institute for Dynamics of Complex Technical Systems, Sandtorstra{\ss}e 1, D-39106 Magdeburg, Germany {benner@mpi-magdeburg.mpg.de}}, B\"ulent Karas\"ozen\thanks{Institute of Applied Mathematics \& Department of Mathematics, Middle East Technical University, 06800 Ankara, Turkey {bulent@metu.edu.tr}}
}

\date{}

\maketitle

\begin{abstract}

This paper discusses a non-intrusive data-driven model order reduction method that learns low-dimensional dynamical models for
a parametrized shallow water equation. We consider the shallow water equation  in non-traditional form (NTSWE).
We focus on learning low-dimensional models in a non-intrusive way. That means, we assume not to have access to a discretized form of the NTSWE in any form. Instead, we have snapshots that are obtained using a black-box solver.
Consequently, we aim at learning reduced-order models only from the snapshots.
Precisely, a reduced-order model is learnt by solving an appropriate least-squares optimization problem in a low-dimensional subspace.
Furthermore, we discuss computational challenges that particularly arise from the optimization problem being ill-conditioned.
Moreover, we extend the non-intrusive model order reduction framework to a parametric case where we make use of the parameter dependency at the level of the partial differential equation.
We illustrate the efficiency of the proposed non-intrusive method to construct reduced-order models for NTSWE and compare it with an intrusive method (proper orthogonal decomposition). We furthermore discuss the predictive capabilities of both models outside the range of the training data.\\


\end{abstract}

\section{Introduction}
Shallow water equations (SWE) are a popular set of hyperbolic PDEs with the capability of describing geophysical wave phenomena, e.g., the Kelvin and Rossby waves in the atmosphere and in the oceans. They are frequently used in geophysical flow prediction \cite{cotter2012mixed}, investigation of baroclinic instability \cite{boss1996stability,vallis2017atmospheric}, and planetary flows \cite{warneford2014thermal}. In this paper, we study a model order reduction (MOR) technique for SWE. MOR techniques allow us to construct low-dimensional models or reduced-order models (ROMs) for a large-scale dynamical system. We refer to the books \cite{Benner17,Rozza14} for an overview of the available techniques. These ROMs are computationally efficient and accurate, and are worthy when a full order model (FOM) needs to be simulated multiple-times for different parameter settings. Additionally,  ROMs are even more valuable in the case of SWE,
when the interest lies
in simulating the model for a very long time horizon.
MOR problems for SWE have been intensively studied in the literature, see, e.g., \cite{Bistrian15,Bistrian17,Esfahanian09,Lozovskiy17,Lozovskiy16,cstefuanescu2014comparison}.

Most MOR techniques are intrusive in nature. This means that these methods require access to the large-scale semi-discretized FOM, preferably in a matrix-vector form. Moreover, ROMs are typically constructed by projecting the high-fidelity FOM onto a low-dimensional subspace using appropriate projection matrices. The proper orthogonal decomposition (POD) is arguably one of the most popular methods that can be seen as a data-driven intrusive method. POD is data-driven in the sense that we require training data that are usually the solution trajectories for given inputs, initial conditions, and parameters. By taking the singular value decomposition (SVD) of the training data, we determine a low-dimensional subspace,  where the most important system dynamics reside.  The efficiency of the ROM is based on the separation of the offline cost for evaluating the FOMs and the online cost for evaluating the ROMs.

One of the major drawbacks of intrusive methods is that they require access to the FOM. However, for a complex dynamical process, it is a challenging task to obtain an explicit discretized FOM.  It is even intractable if the process is simulated using proprietary software.  Therefore, in this work, we are interested in a non-intrusive approach to construct ROMs, where we do not have access to a discretized FOM. We rather have only simulation data, potentially obtained using proprietary software, corresponding to the FOM.

Building a model using only the simulation data directly fits the philosophy of machine learning and neural networks. Using neural networks, a large class of functions \cite{Hornik89} can be approximated. These methods aim at constructing an input-output mapping based on data. They learn a model based on the training data that neither requires explicit access to the high-fidelity model operators nor any additional information about the process. However, the amount of data required to learn the model accurately imposes a burden in the context of large-scale PDE simulations \cite{swischuk2019projection}. Moreover, some ideas from compressive sensing have been used to learn the operators of a FOM from a large library of candidate functions \cite{Rudy19}. However, the success of the method heavily depends on the built library, and we generally need to perform computations in the full-order system dimension, thus making the method very challenging in large-scale settings. In recent years, the operator inference (OpInf) framework to construct ROMs has gained much attention. The framework utilizes the knowledge of nonlinear terms at the  PDE level. In this framework, the operators defining the ROM can be learnt by formulating an optimization problem, without necessitating the discretized operators of the PDEs. Such a scheme was first investigated in \cite{Pehersdorfer16} for polynomial nonlinearities. The methodology was later extended to a class of nonlinear systems that can be written as a polynomial or quadratic-bilinear (QB) system by introducing new state variables in \cite{Kramer19a, QKPW2019_lift_and_learn}. Recently, the authors in \cite{BGKPW2020_OpInf_nonpoly} have extended the approach to nonlinear systems in which the structure of the nonlinearities is preserved while learning ROMs from data.

In this work, we discuss an application of the OpInf framework \cite{Pehersdorfer16} to the parameterized NTSWE. OpInf is also investigated in \cite{Pehersdorfer16} for parametric cases, where the ROMs are constructed at each training parameter via interpolation. However, in this work, we discuss an OpInf framework for the parametric case, where we make use of the known parametric dependency at the PDE level. In the case of a large amount of data, the optimization problem that yields the reduced operators is generally a discrete ill-posed least-squares problem. To mitigate this issue, a regularized least-squares optimization problem is proposed in \cite{Pehersdorfer16}. In this paper, we discuss alternative approaches such as Tikhonov regularization, truncated SVD and truncated QR.

The remaining structure of the paper is as follows. In Section 2, the NTSWE is briefly described. In Section~3, we discuss the OpInf method to infer reduced operators from data and present its extension to the parametric case. Furthermore, we investigate computational issues related to the optimization problem that learns the reduced operators. In Section 4, we present numerical experiments, where ROMs for (parametric) NTSWEs models are inferred directly from data. The inferred ROMs are compared with ROMs obtained from the intrusive POD method. We show that non-intrusive ROMs  outperform in most instances, particularly, in the prediction outside the training data. In Section 5, we provide concluding remarks.

\section{Shallow Water Equation}
\begin{figure}[htb]
	\centering
	\includegraphics[width=0.34\linewidth]{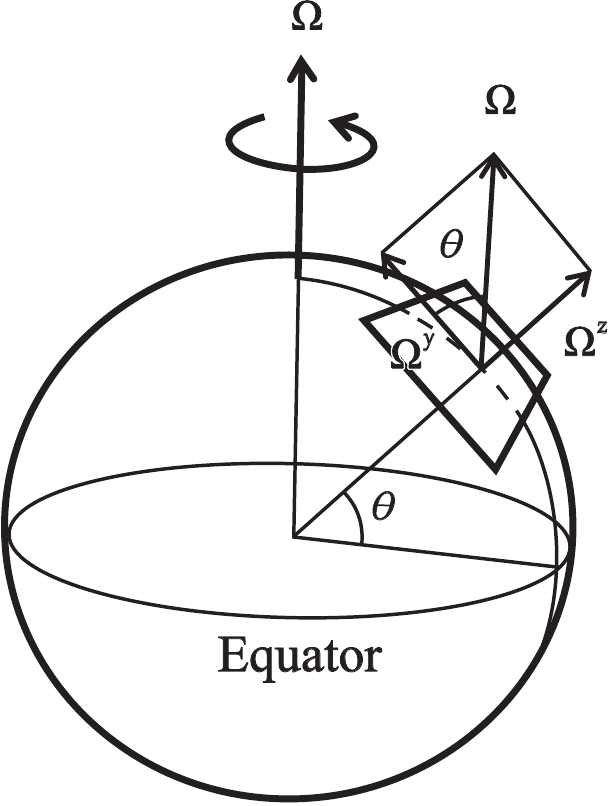}
	\caption{Outline of the components of the angular velocity vector $\Omega$ with respect to the latitude angle $\theta$ as $\Omega^x=0$.}
	\label{fig:cori}
\end{figure}

% Most of the models of the ocean and atmosphere primarily include only the contribution to the Coriolis force from the component of the planetary rotation vector that is locally normal to geopotential surfaces, which is "traditional" approximation. However planetary scale fluid ow in the oceans and in the atmosphere is typically dominated by the Coriolis force. In the following, we describe briefly the non-traditional shallow water equation(NTSWE) \cite{Stewart16,Dellar05,Stewart10}  that includes the full Coriolis force and topography.
In most  ocean and atmosphere models, the Coriolis force only depends  on the component of the planetary rotation vector. It is locally normal to the geopotential surfaces, which is called traditional approximation (TA)~\cite{eckart1960}. The TA is applicable when the horizontal length scales of rotational geophysical flows are much larger than the vertical length scales \cite{gerkema2008}. However,  many atmospheric and oceanographic phenomena are substantially influenced by the non-traditional components of the Coriolis force \cite{stewart2013}, such as deep convection~\cite{marshall99}, Ekman spirals~\cite{leibovich1985}, and internal waves~\cite{gerkema2005}. The NTSWE \cite{Dellar05, Stewart10, Stewart16}  includes the non-traditional components of the Coriolis force with the bottom topography. The non-dimensional NTSWE is governed by the following PDE system  \cite{Stewart16}:
\begin{subequations}\label{eq:ntswe}
\begin{align}
\frac{\partial \mathbf{\tilde{u}}}{\partial t}+qh\mathbf{z}\times\mathbf{u}+\nabla\Phi&=0, \label{eq:ntswe_1}\\
\frac{\partial h}{\partial t}+\nabla\cdot(h\mathbf{u})&=0,
\end{align}
\end{subequations}
where $ \mathbf{\tilde{u}}=:(\tilde{u},\tilde{v}) $ is the canonical velocity,  $ \mathbf{u}=:(u, v) $ is the particle velocity, $ h $ is the height field, $ q $ is the potential vorticity defined as  $q=:\frac{1}{h}\left(\Omega^{z}+\tilde{v}_x-\tilde{u}_y\right)$, and $ \Phi $ is the Bernoulli potential, given by
\begin{align}
\Phi&=\frac{1}{2}\left(u^2+v^2\right)+h_b+h+\frac{1}{2}\delta h\left(\Omega^{x}v-\Omega^{y}u\right).
\end{align}
The particle velocities are given in terms of the canonical velocities  as
%\begin{subequations}
\begin{align}\label{eq:canon}
u& =\tilde{u} - \delta\Omega^y\left(h_b+\frac{1}{2}h\right),\qquad  v = \tilde{v} + \delta \Omega^x\left(h_b+\frac{1}{2}h\right),
\end{align}
%\end{subequations}
where $ h_b$ is the bottom topography, $\delta:=L/R_d$ is the so-called non-traditional parameter with $L$ being the layer thickness scale, and $ R_d $ is the Rossby deformation radius \cite{Stewart10,Dellar05}. $\Omega^x$ and $\Omega^y$ are the $x$ and $ y $ components of the angular velocity vector $\Omega$, respectively, and $x$ and $y$ denote horizontal distances within a constant geopotential surface. The orientation of the $x$ and $y$ axes are considered arbitrary with respect to the North.  %IThe TA of SWE corresponds to setting $ \delta=0 $.
%The effective gravity $g$ is acting in the $z$ direction, perpendicular to the geopotential surface.
The NTSWE \eqref{eq:ntswe} describes inviscid fluid that flows over the bottom topography at $z=h_b(x, y)$ in a frame rotating with angular velocity vector
 ${\mathbf \Omega} = \left(\Omega^{x}, \Omega^{y}, \Omega^{z}\right)$. Both $\Omega^{x}$ and  $\Omega^{y}$ depend on $x$ and $y$ axes but not on $z$. The dimensionless angular velocity vector can be given as~\cite{Stewart12}:
\begin{align} \label{angvel}
\Omega^x=\cos\theta \sin \phi, \quad \Omega^y=\cos\theta \cos \phi, \quad \text{and}\quad \Omega^z=\sin \theta,
\end{align}
where $ \theta $ is the angle corresponding to  the latitude, and $\phi $ is the angle determining the orientation between the $x$-axis and the eastward direction. In this paper, we set the $ x $-axis of the rotation vector to zero, implying that it is aligned to the East. Moreover, we consider the layer thickness scale for the ocean: $ L=1000 $\texttt{m} , the deformation radius $ R_d \approx 6.88 $\texttt{km}, with no bottom topography $h_b=0 $ and the non-traditional parameter $ \delta=0.145 $. In \figurename{\ref{fig:cori}}, the components of the angular velocity vector are shown for the latitude $ \theta $ and $\phi=0$.

\section{Learning Parameterized ROMs of Shallow Water Equations}
In this section, we study an approach to learn ROMs for NTSWE from data, e.g., obtained from proprietary software, or real-world measurements. In Subsection \ref{subsec:QB_SWE}, we begin our discussion with the quadratic form of the parametric NTSWE. We furthermore discuss the construction of ROMs via an intrusive POD method.  Subsection~\ref{subsec:OpInf_para_SWE} presents an operator inference approach to learn the reduced parametric operators  from (simulation) data, where we make use of the knowledge of the parametric dependency at the PDE level. Moreover, we discuss computational aspects for constructing the reduced operators in Subsection~\ref{subsec:ComputationalAspect}.

\subsection{Exploiting the quadratic form of the parameterized NTSWE}\label{subsec:QB_SWE}
The success  of the OpInf approach \cite{Pehersdorfer16,QKPW2019_lift_and_learn} lies in exploiting the structure at the PDE level. The OpInf approach aims at determining a ROM, without having access to the FOMs in a matrix-vector form. It consists of  setting up an optimization problem for determining the reduced operators by taking advantage of the underlying structure of the FOM but at the PDE level.

The NTSWE \eqref{eq:ntswe} can be explicitly written  in terms of the canonical velocities  by taking into account that the $x$-axis of the angular velocity vector is aligned to the East and setting $\phi=0$ as

\begin{subequations}\label{eq:exp_fom}
\begin{align}
\frac{\partial }{\partial t}
\tilde{u}&=-h_x+\sin\theta \ \tilde{v}- \tilde{u}\tilde{u}_x- \tilde{v}\tilde{u}_y + \delta\cos\theta (h\tilde{u})_x- \frac{3}{8}\left(\delta\cos\theta\right)^2(h^2)_x,  \\
\frac{\partial }{\partial t}
\tilde{v}&=-h_y+\sin\theta \ \tilde{u}+\frac{1}{2}\delta\sin\theta\cos\theta \ h- \tilde{u}\tilde{v}_x- \tilde{v}\tilde{v}_y + \delta\cos\theta \left((h\tilde{u})_y+\frac{1}{2}h \left(\tilde{v}_x-\tilde{u}_y\right)\right)- \frac{3}{8}\left(\delta\cos\theta\right)^2(h^2)_y,\\
\frac{\partial }{\partial t}
h&= - (h \tilde{u})_x - (h \tilde{v})_y+\frac{1}{2}\delta\cos\theta( h^2)_x.
\end{align}
\end{subequations}
\subsection{Operator inference approach to learning parameterized reduced-operators}\label{subsec:OpInf_para_SWE}
Let us consider a parameter vector $\mu \in \mathcal{D}\subset \mathbb{R}^{d}$, the state vector $ \mathbf{w}:[0,T]\times \mathcal{D}\rightarrow \mathbb{R}^{N}$ with $N$ degrees of freedom, and the time $ t \in [0,T] \subset \mathbb{R} $. The NTSWE model \eqref{eq:exp_fom} includes linear   terms and quadratic polynomial nonlinearities that can be exploited to create quadratic ROMs. Hence, we consider the following linear-quadratic ODE system:

\begin{align}\label{fom}
\dot{\mathbf{w}}(t;\mu)= \bA(\mu)\mathbf{w}(t;\mu)+\bH(\mu)(\mathbf{w}(t;\mu)\otimes \mathbf{w}(t;\mu)),
\end{align}
where $\bA(\mu)\in \mathbb{R}^{N\times N}$ corresponds to the linear terms, $ \bH(\mu)\in \mathbb{R}^{N\times N^2} $ to the quadratic term. We allow the initial condition to depend on the parameter $\mu$ as well, i.e., $ \mathbf{w}(0,\mu)=\mathbf{w}_0(\mu) $.

Our primary goal is to construct a reduced-order parametrized model that captures the important dynamics of the high-dimensional model \eqref{fom} for a given parameter range as follows:
\begin{align}\label{eq:intrusive}
\dot{\widehat{\mathbf{w}}}(t;\mu)=\widehat{\bA}(\mu)\widehat{\mathbf{w}}(t;\mu)+\widehat{\bH}(\mu)(\widehat{\mathbf{w}}(t;\mu)\otimes\widehat{ \mathbf{w}}(t;\mu)),
\end{align}
where $ \widehat{\bA}(\mu)\in \mathbb{R}^{r\times r}$ and  $\widehat{\bH}(\mu)\in \mathbb{R}^{r\times r^2}$ with $r \ll N$. If the FOM is available in explicit form, e.g., in matrix-vector form,  then intrusive MOR techniques can be applied, such as POD \cite{benner2015two,Benner18} and interpolation-based methods \cite{morBenG19}. Assuming an explicit form of the FOM is available, the ROM can be  constructed with the  projection matrix $\bV \in \R^{N\times r}$ so that $\bw(t,\mu)\approx \bV \widehat\bw(t,\mu),\text{ for all } t \geq 0$ and $\mu \in \cD$  obtained by POD. Then, reduced-order matrices of the system \eqref{eq:intrusive} can be computed as follows:
\begin{equation}\label{eq:intr_romMatx}
 \widehat{\bA}(\mu)=\bV^T\bA(\mu)\bV \in \mathbb{R}^{r\times r},\quad \widehat{\bH}(\mu)=\bV^T\bH(\mu)(\bV\otimes \bV).
\end{equation}
Furthermore, assume that the system matrices in \eqref{fom} depend affinely on functions of the parameter $\mu$:
\begin{subequations}\label{eq:afficeMatx}
  \begin{align}
    \bA(\mu) &= \alpha_1(\mu)\bA_1 + \dots +\alpha_{\mathbf{n_a}}(\mu) \bA_{\mathbf{n_a}},\\
    \bH(\mu) &= \eta_1(\mu)\bH_1 + \dots +\eta_{\mathbf{n_h}}(\mu) \bH_{\mathbf{n_h}},
%    \bB(\mu) &= \beta_1(\mu)\bB_1 + \dots +\beta_{\mathbf{n_b}}(\mu) \bB_{\mathbf{n_b}},
  \end{align}
\end{subequations}
where $\bA_i \in \R^{N\times N}$, $\bH_j \in \R^{N\times N^2}$ are constant matrices, and $\alpha_i(\mu),\eta_j(\mu): \R^d \rightarrow \R$ are smooth functions of the parameter $\mu$.    In this case, the reduced-matrices in \eqref{eq:intr_romMatx} can be precomputed, e.g., $\widehat{\bA}(\mu)  = \alpha_1(\mu)\widehat{\bA}_1 + \dots +\alpha_{\mathbf{n_a}}(\mu) \widehat{\bA}_{\mathbf{n_a}}$, where $\widehat\bA_i = \bV^T\bA_i\bV \ , \ i \in \{1,\ldots,\mathbf{n_a}\}$.

However, as discussed earlier, it is not easy or almost impossible to obtain the FOM in an explicit matrix form, from  proprietary software. Therefore, our primary interest lies in constructing reduced-order operators without having access to the FOM, but rather having access only to simulation data and some knowledge at the PDE level. With this aim, we collect simulation data for a training parameter set, $ \mu_i \in \mathcal{D} $ for $ i=1,\dots,M$. Thus, let us define a global snapshot matrix:
\begin{align}\label{eq:global_snap}
\bS_{\mu}=\left[\bS(\mu_1),\dots,\bS(\mu_M)\right], \quad
\bS(\mu_i)=\left[\mathbf{w}(t_1;\mu_i),\mathbf{w}(t_2;\mu_i),\ldots,\mathbf{w}(t_K;\mu_i)\right]\in \mathbb{R}^{N\times K},
\end{align}
where $\bw(t_j,\mu_i)$ denotes the value at time $t_j$ for the parameter $\mu_i$. The projection matrix $\bV$ is determined by the SVD of the snapshot matrix
\begin{equation}
 \bS_\mu = \bV_\mu\Sigma_\mu\bU_\mu^T,
\end{equation}
where $\bV_\mu\in \R^{N\times M\cdot K},\Sigma_\mu\in \R^{M\cdot K\times M\cdot K}, \bU_\mu \in \R^{M\cdot K\times M\cdot K}$,  and  $\bV$ is given  then by the first $r$ columns of $\bV_\mu$.
In order to determine reduced operators by employing an OpInf approach, we first project the snapshot matrix $\bS_\mu$ onto the dominant subspace spanned by $\bV$, yielding the reduced snapshot matrix:

\begin{equation}
 \widehat\bS_\mu := \bV^T \bS_\mu = \begin{bmatrix}\widehat\bS(\mu_1),\dots,\widehat\bS(\mu_M)\end{bmatrix},
\end{equation}
where
\begin{align*}
\widehat\bS(\mu_i)=\left[\widehat\bw(t_1;\mu_i),\widehat\bw(t_2;\mu_i),\ldots,\widehat\bw(t_K;\mu_i)\right]\in \mathbb{R}^{r\times K}
\end{align*}
in which $\widehat\bw(t_j,\mu_i) := \bV^T\bw(t_j,\mu_i)$. Furthermore, let us define
\begin{align}
\dot{\widehat{\bS}}_\mu =  \begin{bmatrix}\dot{\widehat\bS}(\mu_1),\dots,\dot{\widehat\bS}(\mu_M)\end{bmatrix},
\end{align}
where $\dot{\widehat\bS}(\mu_i)$ can either be determined using the right-hand side of \eqref{fom}--if accessible--followed by projecting using $\bV$, or can be approximated using $\widehat{\bS}(\mu_i)$ by employing a time-derivative approximation scheme, see, e.g.,~\cite{Pehersdorfer16}. Subsequently, the reduced operators of the reduced parametric model \eqref{eq:intrusive} are determined by solving the following least-squares problem:
\begin{equation}\label{eq:opt_opeInf}
 \min_{\widehat\bA_i\in \Rrr,\widehat\bH_j\in\R^{r\times r^2}} \sum_{k = 1}^M \left\|  - \dot{\widehat{\bS}}\left(\mu_k\right)^T + \sum_{i = 1}^{\mathbf{n_a}}\left(\alpha_i(\mu_k)\widehat\bS(\mu_k)^T\widehat\bA^T_i\right) +  \sum_{i = 1}^{\mathbf{n_h}}\left(\eta_i(\mu_k)\left(\widehat\bS(\mu_k)\hat{\otimes}\widehat\bS(\mu_k)\right)^T\widehat\bH^T_i\right)  \right\|_{F}^{2}
\end{equation}
where $\hat{\otimes}  $ denotes the column-wise Kronecker product.
 It can be noted that the optimization problem \eqref{eq:opt_opeInf} does not involve any explicit knowledge  of the FOM, it only involves   simulation data projected onto the dominant POD subspace. Moreover, we can rewrite the optimization problem \eqref{eq:opt_opeInf} in standard form as follows:
\begin{equation}\label{eq:lsqr}
 \min_{\cX \in \R^{r\times \mathbf{n_a} r+ \mathbf{n_h}r^2}}\sum_{k=1}^M\left\|\cA(\mu_k)\cX^T - \dot{\widehat{\bS}}\left(\mu_k\right)^T\right\|_{F}^{2},
\end{equation}
where
\begin{align*}
 \cX &= \begin{bmatrix} \widehat\bA_1,\ldots,\widehat\bA_\mathbf{n_a},\widehat\bH_1,\ldots,\widehat\bH_\mathbf{n_h} \end{bmatrix},~\text{and}\\
 \cA(\mu_k) &= \begin{bmatrix} \left[\alpha_1(\mu_k),\ldots, \alpha_{\mathbf{n_a}}(\mu_k)\right]\otimes\widehat\bS(\mu_k)^T,\left[\eta_1(\mu_k),
 \ldots,\eta_{\mathbf{n_h}}\right]\otimes\left(\widehat\bS(\mu_k)\hat{\otimes}\widehat\bS(\mu_k)\right)^T
\end{bmatrix}.
\end{align*}
%In the following, we discuss computational aspects of the optimization problem  \eqref{eq:lsqr}.

\subsection{Computational Aspects}\label{subsec:ComputationalAspect}

In this section, we discuss computational aspects of the OpInf approach   \eqref{eq:lsqr}.
Solving the least-squares problem \eqref{eq:lsqr} can be a computationally challenging task because  (a) the problem can be highly ill-conditioned, and (b) its computational cost grows quadratically with the order $r$ of the reduced system and linearly with the number of snapshots. The computational cost of the optimization problem \eqref{eq:lsqr} can be reduced by decoupling of the least-squares problem. In the remainder, we discuss techniques for the conditioning.
In the OpInf framework, inferred operators are solutions of the potentially  discrete ill-posed least-squares problem \eqref{eq:lsqr}, where ill-conditioning may be due to nearly linearly dependent columns of the snapshot matrix.
When the distance to a matrix with linearly dependent columns decreases, the condition number of the data matrix in the least-squares problem increases. Hence, the least-squares problem arising from the OpInf framework needs a suitable regularization method.
There exist different ways to deal with this issue.

 A suitable and widely used candidate for this task is Tikhonov regularization \cite{tikhonov77}. Tikhonov regularization filters the small singular values to reduce the amplification effect on the least-squares algorithm. The quality of learning via Tikhonov regularization depends on the L-curve \cite{hansen2001curve}. Using the L-curve information, the stability of the learning algorithm can be improved \cite{kramer19b}.  However, the computation of the efficient Tikhonov parameter is costly for large problems \cite{calvetti2002curve}. Although one can argue that the problem \eqref{eq:lsqr} is in a low dimension, it still can be of a large scale when the number of snapshots or/and the number of training parameters are large.
The  Tikhonov regularization applied to \eqref{eq:lsqr}  can be written in a compact form as follows:
\begin{equation}\label{eq:tikhonov}
\min_{{x_i} \in \mathbb{R}^{\mathbf{n_a}r+ \mathbf{n_h}r^2}}\left\|\mathcal{A}_\mu x_i - s_i\right\|^2_{2}+\lambda\left\|  x_i \right\|_2^2, \quad i=1,\ldots,r,
\end{equation}
where $ x_i $ are the columns of $ \mathcal{X}^T $, $ s_i $ are the columns of $ \dot{\widehat{\mathbf{S}}}_\mu $ and $ \mathcal{A}_\mu=\left[\mathcal{A}^T(\mu_1),\dots,\mathcal{A}^T(\mu_M)\right]^T $.
A heuristic approach to deal with the conditioning of the data matrix is proposed in \cite{Pehersdorfer16} in which a subset of the data is considered,  by taking the data in a regular interval, e.g., every $10$th time-step. In \cite{Pehersdorfer16},  it is shown that this can alleviate the ill-conditioning problem to some extend  in some cases. However, the choice of the interval should be done in such a way that the important snapshots are not missed, thus the choice of the interval plays a key role. This problem can be referred to as a heuristic column subset selection problem (CSSP). The CSSP seeks to find a subset of the most linearly independent columns of a matrix which gives the best information in the matrix.

The truncated QR method (tQR) with the minimum norm solution \cite{chan1992some} finds a suitable subset for the CSSP, which can be used to find an accurate solution of the rank-deficient least-squares problems. For the tQR method, we use the QR decomposition of the data matrix with column pivoting (QR-CP). A major advantage of the tQR algorithm is that it allows us to monitor the linearly dependent columns via QR-CP. Thus, it also allows us to improve the condition number of the data matrix by selecting linearly independent vectors from the data matrix. The QR-CP algorithm naturally finds a subset for the CSSP problem via a permutation matrix. One alternative to this approach is the truncated SVD (tSVD) algorithm. The tSVD is also an attractive method with its best rank-$k$ approximation. The tSVD method and tQR method generally give very close solutions. Nevertheless, the calculation of the  tSVD is more expensive than the tQR algorithm.
%The QR-CP algorithm generally results in a  reliable a low-rank approximation with a few rare exceptions such as the Kahan matrix \cite{golub2013matrix}. Though column pivoting adds an additional computational cost to the QR factorization algorithm, the computational cost of QR-CP \emph{almost} becomes as cheap as the QR method with the help of new randomized QR-CP methods \cite{duersch2017randomized}.

For simplicity, suppose $ \cA_\mu \in \mathbb{R}^{m\times n} $ is exactly rank deficient with $ \text{rank}(\cA_\mu) = p $. Then, there always exists a QR-CP factorization of $ \cA_\mu $ of the form
\begin{align}\label{eq:QR}
\cA_\mu \Pi=QR,
\end{align}
where $ \Pi\in \mathbb{R}^{n\times n}  $ is a permutation matrix, $ Q \in \mathbb{R}^{m\times m}  $  is an orthogonal matrix, and $ R $ is an upper triangular matrix of the form
\begin{align}\label{eq:R}
R=\begin{pmatrix}
R_1&R_2\\
0&0
\end{pmatrix}
\end{align}
and $ R_1 \in \mathbb{R}^{p\times p} $ is upper triangular with $ \text{rank}(R_1)=p $. The diagonal entries of $ R $ in \eqref{eq:QR} satisfy $ |R_{ii}| \geq |R_{jj}|$ with $j>i$ so that the effective rank of $ \cA $ can be determined as the smallest integer $ p $ such that
\begin{align*}
|R_{p+1,p+1}|<\texttt{tol}\cdot \ |R_{11}|,
\end{align*}
where \texttt{tol} can be considered as the tolerance for the linear dependency of the columns of data matrix $ \cA_\mu $. Hence, this approach gives us a flexibility of adjustment on the condition number of the data matrix $ \cA_\mu $. The tolerance of tQR $ \texttt{tol} $ can also be determined by the L-curve \cite{Roeck02}. Among all solutions of the optimization problem \eqref{eq:lsqr}, we take the unique minimum $L_2$ norm solution. The minimum norm solution of \eqref{eq:lsqr} by tQR method can be obtained as in \cite{chan1992some}. In this paper, we study an QR-CP based regularization and compare with the Tikhonov regularization in the next section.

\section{Numerical results}
In this section, we demonstrate the performance of the OpInf approach for two numerical test problems and compare it with the intrusive POD method.
%It is initially in geostrophic balance which particularly means that the pressure gradient balanced by the traditional component of the Coriolis force '$ \Omega^z $'.
We also study the prediction capabilities of the OpInf approach for the parametric and non-parametric NTSWE. Furthermore, we examine the non-intrusive methods with the Tikhonov regularization \eqref{eq:tikhonov} for the penalty parameter $ \lambda=0.01 $ and tQR method with the tolerance $10^{-6} $. To determine the regularization parameters of Tikhonov regularization and tQR, we used L-curve criteria. As a first test example, we consider the propagation of the inertia-gravity waves by Coriolis force, known as  geostrophic adjustment \cite{Stewart16}. In the second example, we investigate the shear instability in the form of a roll-up of an unstable shear layer, known as barotropic instability  \cite{Stewart16}.

The NTSWE \eqref{eq:exp_fom} is semi-discretized in space by replacing the first-order spatial derivatives with central finite-differences. The resulting system is a quadratic semi-discrete system that depends on the parameter $\mu=\theta$:
\begin{align}\label{discrete_fom}
\dot{\mathbf{w}}(t;\mu)= \bA(\mu)\mathbf{w}(t;\mu)+\bH(\mu)(\mathbf{w}(t;\mu)\otimes \mathbf{w}(t;\mu)),
\end{align}
where $\bA(\mu)\in \mathbb{R}^{N\times N}$ corresponds to the linear terms, $ \bH(\mu)\in \mathbb{R}^{N\times N^2} $ to the quadratic term, $\mu \in \mathcal{D}\subset \mathbb{R}^{d}$, the state vector $ \mathbf{w}:[0,T]\times \mathcal{D}\rightarrow \mathbb{R}^{N}$ with $N$ degrees of freedom, and the time $ t \in [0,T] \subset \mathbb{R} $. The operators $ \bA(\mu) $  and $ \bH(\mu) $ in \eqref{eq:afficeMatx} have affine parameter dependence \eqref{eq:afficeMatx} with respect to the parameter $ \mu $ as follows:
\begin{subequations}\label{eq:paramFunc}
	\begin{align}
	&\alpha_1(\mu)=1,\quad \alpha_2(\mu)=\Omega^z=\sin(\mu),\quad\alpha_3(\mu)=\Omega^z\Omega^y=\sin(\mu)\cos(\mu),\\
	&\eta_1(\mu)=1,\quad \eta_2(\mu)=\Omega^y=\cos(\mu),\quad\eta_3(\mu)=(\Omega^y)^2=(\cos(\mu))^2.
	\end{align}
\end{subequations}

We consider the NTSWE under periodic boundary conditions, and assume there are no forcing or input terms in \eqref{discrete_fom}. All the models are simulated using the function \texttt{ode15s} in \matlab~with both the relative and absolute error tolerances set to $10^{-8} $. In all numerical examples, the spatial domain is discretized with $ 101\times101 $ equidistant grid points. The snapshots are sampled at equidistant time instances with the time step $\Delta t =0.1$.
The accuracy of the ROMs is measured using the relative error in the Frobenius norm
\begin{align}\label{eq:rel_err}
\mathcal{E}=\frac{||\bS_{\text{FOM}}-\bV\bS_{\text{ROM}}||_{F}}{||\bS_{\text{FOM}}||_{F}},
\end{align}
where $ \bS_{\text{FOM}}\in\mathbb{R}^{N\times K} $ is the snapshot matrix of the FOM and  $\bS_{\text{ROM}}\in\mathbb{R}^{r\times K} $ is the snapshot matrix of either the non-intrusive or the intrusive ROMs. In the parametric case, the relative error is computed with snapshot matrices by concatenating the trajectories for parameter samples. A typical approach to determine the reduced dimension $ r $ is done through the projection error as follows:
\begin{align}\label{eq:proj_err}
\mathcal{E}_{\text{proj}}=\frac{||\bS_{\text{FOM}}-\bV\bV^T\bS_{\text{FOM}}||_{F}}{||\bS_{\text{FOM}}||_{F}}.
\end{align}

\subsection{Single-layer geostrophic adjustment}
The initial conditions are prescribed in the form of a motionless layer with an upward bulge of the height field in a periodic domain $ [-5, 5]\times [-5,5]$ :
\begin{align*}
&h(x, y, 0) = 1+\frac{1}{2} \exp\left[-\left(\frac{4x}{5}\right)^2-\left(\frac{4y}{5}\right)^2\right], \qquad
u(x, y, 0) = 0,\qquad
v(x, y, 0) = 0.
\end{align*}
The inertia-gravity waves propagate after the collapse of the initial symmetric peak with respect to the axes. Nonlinear interactions create shorter waves, propagating around the domain, and the interactions construct more complicated patterns \cite{Stewart16}.

\FloatBarrier
\subsubsection{Non-parametric case}
We consider the NTSWE for a fixed parameter $\mu = \tfrac{\pi}{4}$. The snapshots corresponding to the FOM are collected by solving \eqref{eq:ntswe} in the time domain $[0, T]$ with $T=60$. The discrete state vectors are concatenated into $\bw=[\bu,\bv,\bh]^T\in\mathbb{R}^{30000}$, leading to a training set of the size $ \bS_{\text{FOM}}\in\mathbb{R}^{30000\times601} $. The inferred non-intrusive ROM and the intrusive ROM are then used to predict the height field outside of the training set at  time $T=80$.
%To show the efficiency of the two regularization methods on the OpInf framework, we considered the reduce dimensions where both methods are stable under %fixed regularization parameter.

\begin{figure}[H]
	\centering
	\setlength\fheight{1.5cm}
	\setlength\fwidth{.5\textwidth}
	{\small
    \tikzsetnextfilename{fig//Nex1/sing}%
    \input{fig//Nex1/sing.tikz}%

	}
	\caption{Single-layer geostrophic adjustment: Normalized singular values.}
	\label{fig:singNex1}
\end{figure}
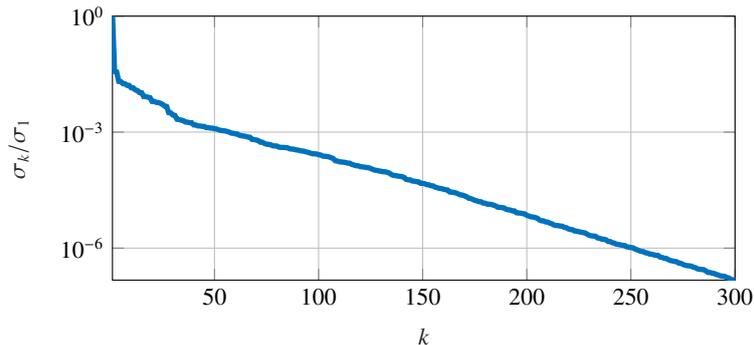

The decay of the leading $300$ normalized singular values of the snapshot matrix $ \bS_{\text{FOM}} $ is shown in  \figurename{ \ref{fig:singNex1}}. The slow decay of the normalized singular values indicates the difficulty of obtaining accurate ROMs for a small number of  POD modes, which is a common problem for hyperbolic PDEs, e.g., like the NTSWE.

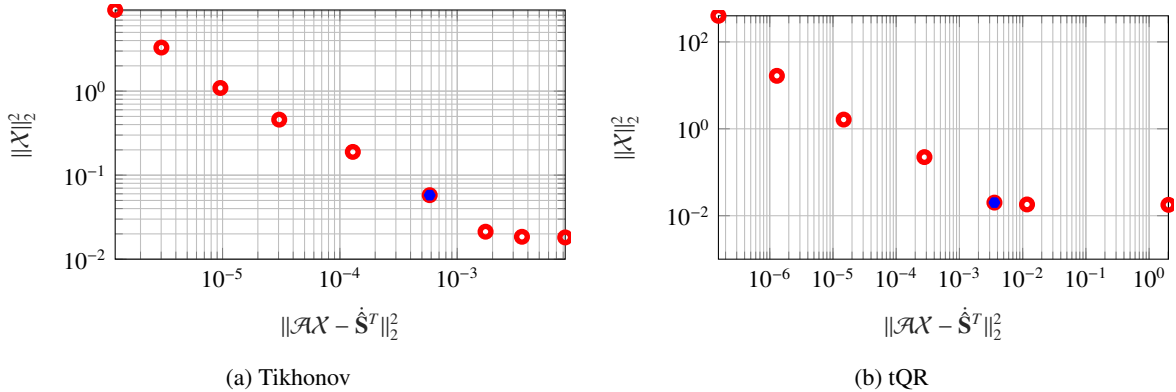
\begin{figure}[H]
		\begin{subfigure}[t]{0.45\textwidth}
		\centering
		\setlength\fheight{1.5cm}
		\setlength\fwidth{.8\textwidth}
		{\small
    \tikzsetnextfilename{fig//tikL}%
    \input{fig//tikL.tikz}%

		}
		
		%\vspace{-0.5cm}
		\caption{Tikhonov}
		\label{fig:tikL}
	\end{subfigure}
	\begin{subfigure}[t]{0.45\textwidth}
		\centering
		\setlength\fheight{1.5cm}
		\setlength\fwidth{.8\textwidth}
		{\small
    \tikzsetnextfilename{fig//tqrL}%
    \input{fig//tqrL.tikz}%

		}
		
		%\vspace{-0.5cm}
		\caption{tQR }
		\label{fig:tqrL}
	\end{subfigure}
\caption{L-curves: Tikhonov (left) , tQR (right).}
\label{fig:Lcurves}
\end{figure}

In \figurename{ \ref{fig:Lcurves}} ,  we demonstrate the L-curves for the basis sizes $r=20$ , tQR tolerances $\texttt{tol}=[10^{-4},10^{-5},\dots,10^{-10}]$ and Tikhonov parameters $\lambda=[10^{1}, 10^{0}, \dots, 10^{-7},]$.  The vertical axis shows the squared norm of the learned operators, and the horizontal axis shows the least-squares residual. The computation of the L-curves for each dimension is costly therefore
the Tikhonov regularization parameter $\lambda$ and the tQR tolerance $ \texttt{tol} $ are chosen close to the upper part of the corners at the L-curves for $r=20$ which produce stable solutions from dimension $r=20$ up to dimension $r=75$. The chosen tQR tolerance and Tikhonov parameter $\lambda $ are shown as blue dot in \figurename{ \ref{fig:Lcurves}}.

%The regularization tolerance $\texttt{tol}$ is chosen as $\lambda$  which lie in the upper part of the corners of the L-curves. Moreover, the chosen tQR tolerance and Tikhonov parameter $\lambda $ are shown as blue dot in \figurename{ \ref{fig:Lcurves}}. For higher dimensions up to $r=75$, we used the same regularization parameters since they are giving meaningful and stable results.

Next, we build a POD-projection matrix $ \bV \in \mathbb{R}^{30000 \times r} $ with $ r=75 $, which leads to a projection error \eqref{eq:proj_err}  $\mathcal{E}_{\text{proj}}= 2.07\cdot 10^{-3}$. Then, we construct ROMs using the intrusive POD and non-intrusive OpInf methods. The non-intrusive methods regularized with both Tikhonov and tQR-based regularizers for comparison. To examine the accuracy of these ROMs, we perform time-domain simulations and compare them with the FOM. In \tablename{ \ref{tab:1}}, the accuracy of the intrusive and non-intrusive ROMs are compared using the FOM-ROM error \eqref{eq:rel_err}. \tablename{ \ref{tab:1}} indicates that the non-intrusive method yields a better ROM as compared to the intrusive POD method, and for this example, both regularizations apparently yield a similar result.

\begin{table}[H]
	\centering
	\begin{tabular}{|c|c|c|c|}
		\hline
		Method & POD & Non-intrusive (Tikhonov regularizer) & Non-intrusive (tQR)\\\hline
$\mathcal{E}$ &$3.27\cdot 10^{-3}$ & $ 2.07 \cdot 10^{-3} $  &  $2.07 \cdot 10^{-3} $ \\
\hline
    \end{tabular}
\caption{Single-layer geostrophic adjustment: Comparison of intrusive and non-intrusive ROMs.}
\label{tab:1}
\end{table}

\begin{figure}[H]
		\centering
		\setlength\fheight{1.5cm}
		\setlength\fwidth{.6\textwidth}
		{\small
    \tikzsetnextfilename{fig//Nex1/rel_err}%
    \input{fig//Nex1/rel_err.tikz}%

		}
	\caption{Single-layer geostrophic adjustment: Relative ROMs errors.}
		\label{fig:rel_errNex1}
\end{figure}
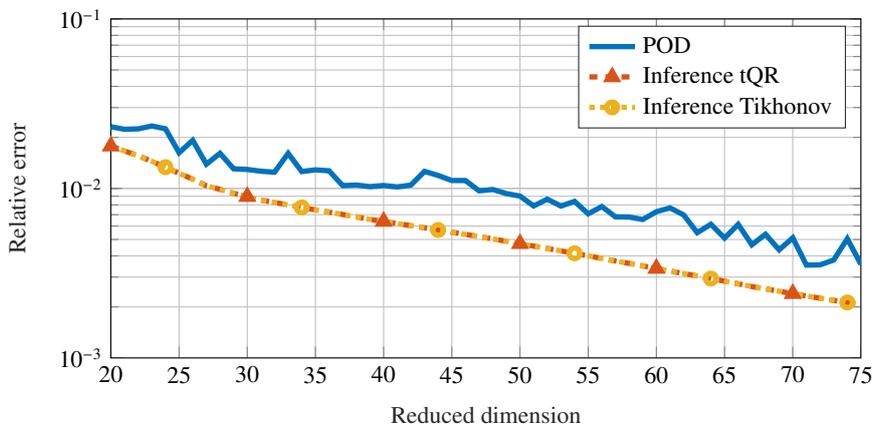

Moreover, we study the quality of the ROMs as the ROM order increases. For this, we compare the relative error $\cE$ \eqref{eq:rel_err} obtained using intrusive and non-intrusive methods in \figurename{ \ref{fig:rel_errNex1}}. We observe that the relative error \eqref{eq:rel_err} using the intrusive ROM does not decrease as smoothly as in the case of the non-intrusive case. Furthermore, we observe that both regularizers perform equally well for all orders of the ROMs.
However, the penalty parameter of the Tikhonov regularization has to be determined by the L-curve, which requires the SVD computation of the data matrix \cite{hansen2001curve}. Thus, the Tikhonov regularization, combined with L-curve information, is more expensive than the tQR.
%Moreover, L-curve information may not give an informative regularization parameter for the least-squares problem \eqref{eq:tikhonov} as discussed in \cite{kramer19b}.

\begin{figure}[H]
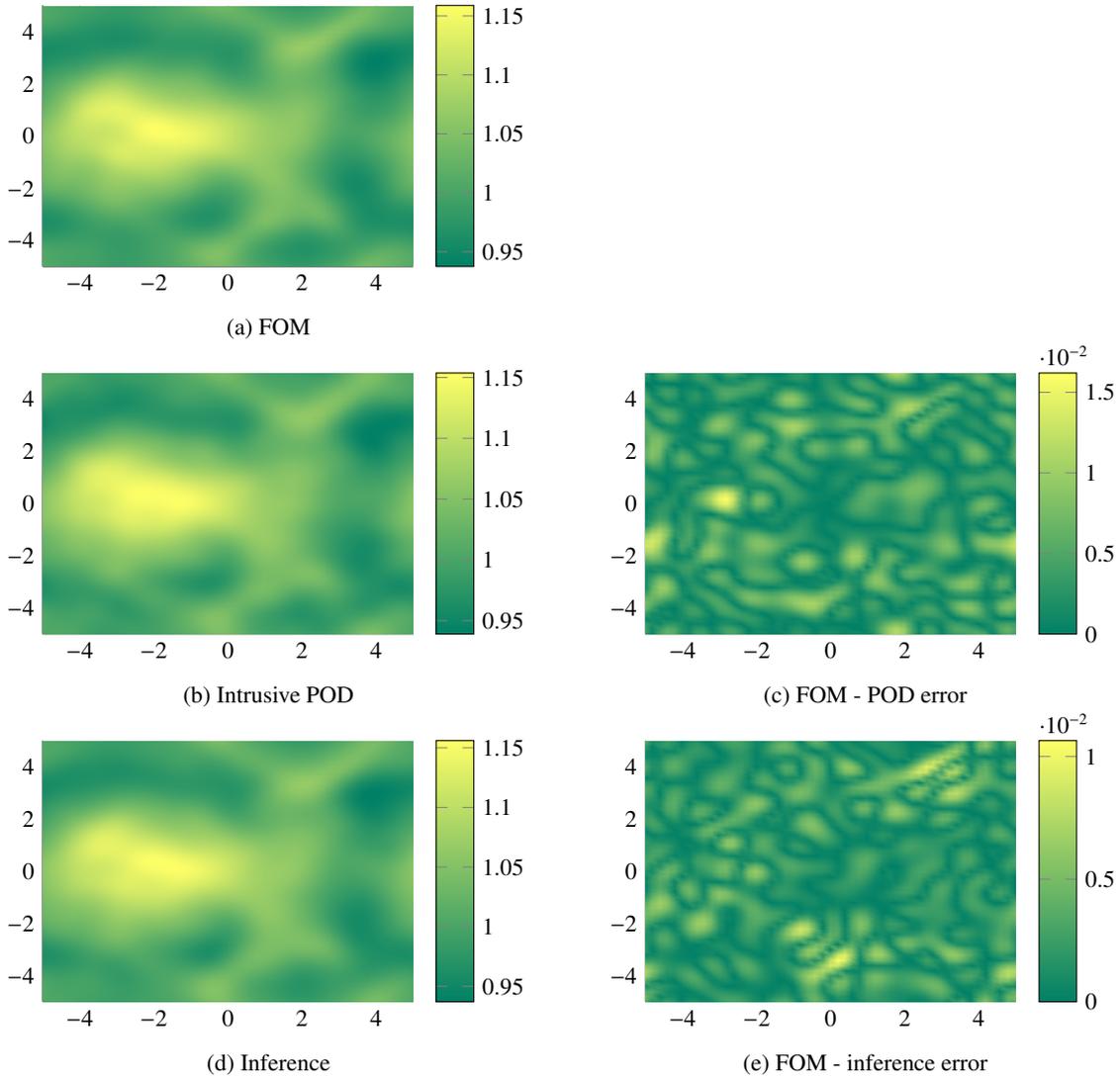

	\begin{subfigure}[t]{0.45\textwidth}
		\centering
		\setlength\fheight{1.5cm}
		\setlength\fwidth{.8\textwidth}
		{\small
    \tikzsetnextfilename{fig//Nex1/fom}%
    \input{fig//Nex1/fom.tikz}%

		}
		\caption{FOM }
		\label{fig:fomNex1}
	\end{subfigure}
\\
	\begin{subfigure}[t]{0.45\textwidth}
		\centering
		\setlength\fheight{1.5cm}
		\setlength\fwidth{.8\textwidth}
		{\small
    \tikzsetnextfilename{fig//Nex1/pod}%
    \input{fig//Nex1/pod.tikz}%

		}
		\caption{Intrusive POD}
		\label{fig:podNex1}
	\end{subfigure}
	\begin{subfigure}[t]{0.45\textwidth}
	\centering
	\setlength\fheight{1.5cm}
	\setlength\fwidth{.7\textwidth}
	{\small
    \tikzsetnextfilename{fig//Nex1/fom_pod}%
    \input{fig//Nex1/fom_pod.tikz}%

	}
	\caption{FOM - POD error}
	\label{fig:fom_podNex1}
	\end{subfigure}
\\
	\begin{subfigure}[t]{0.45\textwidth}
	\centering
	\setlength\fheight{1.5cm}
	\setlength\fwidth{.8\textwidth}
	{\small
    \tikzsetnextfilename{fig//Nex1/inf}%
    \input{fig//Nex1/inf.tikz}%

	}
	\caption{Inference}
	\label{fig:infNex1}
	\end{subfigure}
	\begin{subfigure}[t]{0.45\textwidth}
	\centering
	\setlength\fheight{1.5cm}
	\setlength\fwidth{.8\textwidth}
	{\small
    \tikzsetnextfilename{fig//Nex1/fom_inf}%
    \input{fig//Nex1/fom_inf.tikz}%

	}
	\caption{FOM - inference error}
	\label{fig:fom_infNex1}
	\end{subfigure}
    \caption{Single-layer geostrophic adjustment: Comparison of the height field obtained using the FOM and ROMs of order $r = 75$ at time $ T=60 $.}
    \label{fig:comparion_HF_t60}
\end{figure}

In the rest of this paper, for time-domain simulations and prediction, we show the results for the OpInf with tQR  since both the Tikhonov and tQR-based methods yield comparable results with respect to the projection error \eqref{eq:proj_err}. The height field of the FOM and ROMs of order $r=75$  at time $T=60$ are shown in \figurename{ \ref{fig:comparion_HF_t60}}. These figures show that both ROMs are very close to the FOM solutions, with the non-intrusive ROM being slightly more accurate.

Furthermore, we discuss the prediction capabilities of the obtained ROMs of order $r = 75$. For this, we predict the height field at time $T=80$  in  \figurename{ \ref{fig:prediction_HF_t80}}. Note that we have trained the model using the data up to $T = 60$.
% It shows that non-intrusive ROM performs better than the intrusive ROM.
It indicates that the height field can be predicted with significantly higher accuracy using the non-intrusive ROM as compared to the intrusive POD model.  %We note that the accuracy of the OpInf can be increased by increasing the training interval. %\pawan{Probably it is true for POD model as well.. So please remove of modify it}{}

%\pawan{Please change inference model to non-intrusive model or non-intrusive ROM!!}{}

\begin{figure}[H]
	\begin{subfigure}[t]{0.45\textwidth}
		\centering
		\setlength\fheight{1.5cm}
		\setlength\fwidth{.8\textwidth}
		{\small
    \tikzsetnextfilename{fig//Nex1/Pred/fom}%
    \input{fig//Nex1/Pred/fom.tikz}%

		}
		\caption{FOM}
		\label{fig:fomNex1pred}
	\end{subfigure}
	\\
	\begin{subfigure}[t]{0.45\textwidth}
		\centering
		\setlength\fheight{1.5cm}
		\setlength\fwidth{.8\textwidth}
		{\small
    \tikzsetnextfilename{fig//Nex1/Pred/pod}%
    \input{fig//Nex1/Pred/pod.tikz}%

		}
		\caption{Intrusive POD}
		\label{fig:podNex1pred}
	\end{subfigure}
	\begin{subfigure}[t]{0.45\textwidth}
		\centering
		\setlength\fheight{1.5cm}
		\setlength\fwidth{.8\textwidth}
		{\small
    \tikzsetnextfilename{fig//Nex1/Pred/fom_pod}%
    \input{fig//Nex1/Pred/fom_pod.tikz}%

		}\\
		\caption{FOM - POD error}
		\label{fig:fom_podNex1pred}
	\end{subfigure}
	\\
	\begin{subfigure}[t]{0.45\textwidth}
		\centering
		\setlength\fheight{1.5cm}
		\setlength\fwidth{.8\textwidth}
		{\small
    \tikzsetnextfilename{fig//Nex1/Pred/inf}%
    \input{fig//Nex1/Pred/inf.tikz}%

		}
		\caption{Inference}
		\label{fig:infNex1pred}
	\end{subfigure}
	\begin{subfigure}[t]{0.45\textwidth}
		\centering
		\setlength\fheight{1.5cm}
		\setlength\fwidth{.8\textwidth}
		{\small
    \tikzsetnextfilename{fig//Nex1/Pred/fom_inf}%
    \input{fig//Nex1/Pred/fom_inf.tikz}%

		}
		\caption{FOM - inference error}
		\label{fig:fom_infNex1pred}
	\end{subfigure}
	\caption{Single-layer geostrophic adjustment: Prediction of the height field obtained using the FOM and ROMs of order $r=75$ at time $T=80$.}	
	\label{fig:prediction_HF_t80}
\end{figure}

\FloatBarrier
\subsubsection{Parametric Case}
In this example, we set the parameter domain of NTSWE as $ \mathcal{D} = \left[\frac{\pi}{6},\frac{\pi}{3}\right] \subset \mathbb{R} $. We consider the time domain $ [0, T ] \subset R  $ with final time $ T=10 $. The trajectories are generated with $M=5$ equidistantly distributed parameters
$\mu_1,\mu_2,\ldots,\mu_5 \in \mathcal{D}$. The concatenated snapshot matrix $ \bS_{\text{FOM}} \in\mathbb{R}^{30000\times 505} $ is constructed as in \eqref{eq:global_snap}.

\begin{figure}[H]
	\centering
		\setlength\fheight{1.5cm}
		\setlength\fwidth{.5\textwidth}
		{\small
    \tikzsetnextfilename{fig//Pex1/sing}%
    \input{fig//Pex1/sing.tikz}%

		}
	\caption{Single-layer geostrophic adjustment: Normalized singular values.}
	\label{fig:singPex1}
\end{figure}
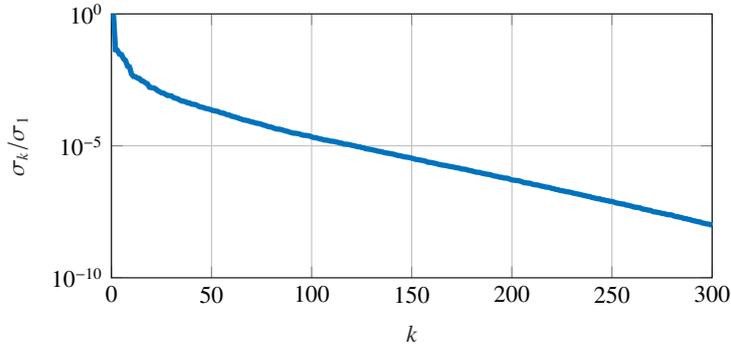

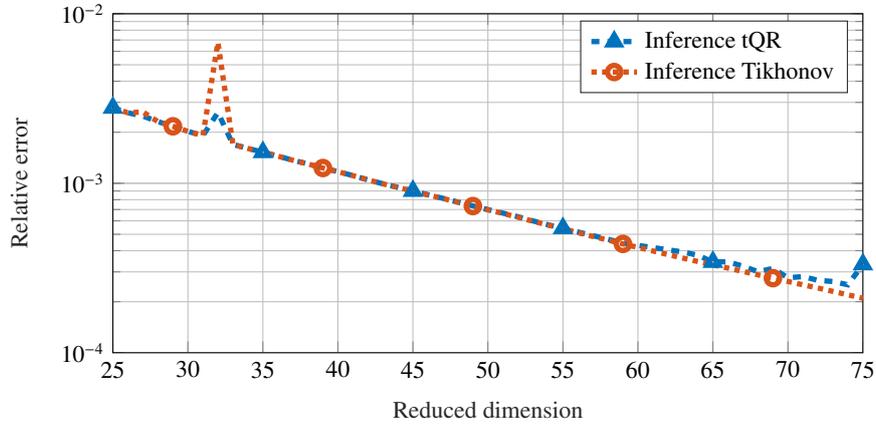
\begin{figure}[H]
		\centering
		\setlength\fheight{1.5cm}
		\setlength\fwidth{.6\textwidth}
		{\small
    \tikzsetnextfilename{fig//Pex1/rel_err_training}%
    \input{fig//Pex1/rel_err_training.tikz}%

		}
		\caption{Single-layer geostrophic adjustment: Relative ROMs errors.}
        \label{fig:rel_err_Pex1}
\end{figure}

The normalized singular values in \figurename{ \ref{fig:singPex1}} decay similar to the non-parametric case in \figurename{ \ref{fig:singNex1}}.
To determine the accuracy of the non-intrusive ROMs, we compute the projection error \eqref{eq:proj_err} in the training set, which is $\mathcal{E}_{\text{proj}}= 2.10\cdot10^{-4} $ for $ r=75 $. In \tablename{ \ref{tab:2}}, the relative errors \eqref{eq:rel_err} are very close to the projection error \eqref{eq:proj_err} for the ROMs of order $ r=75$, which indicates that non-intrusive ROM solutions of both methods have the same level of accuracy. For the parametric case, we present only the results of the non-intrusive approach; however, we observe a similar behavior as in the non-parametric case when non-intrusive and intrusive methods are compared. We also demonstrate the relative errors for the non-intrusive ROMs of orders $r=25$ to $r=75$ over the training set in \figurename{ \ref{fig:rel_err_Pex1}}. Again, we observe the same behavior as for the non-intrusive model; when the order of the ROM increases, the relative errors in the training set decrease.
% For the dimension of $ r=32 $ both regularization schemes provide stable results.

\begin{table}[H]
	\centering
	\begin{tabular}{|c|c|c|}
		\hline
		Method  & Non-intrusive (Tikhonov regularizer) & Non-intrusive (tQR)\\\hline
		$\mathcal{E}$ & $ 2.10 \cdot 10^{-4} $  &  $3.32 \cdot 10^{-4} $ \\
		\hline
	\end{tabular}
	\caption{Single-layer geostrophic adjustment: Comparison of non-intrusive ROMs.}
	\label{tab:2}
\end{table}

\begin{figure}[H]
		\centering
		\setlength\fheight{1.5cm}
		\setlength\fwidth{.7\textwidth}
		{\small
    \tikzsetnextfilename{fig//Pex1/param_err}%
    \input{fig//Pex1/param_err.tikz}%

		}
	\caption{Single-layer geostrophic adjustment: Relative error for testing and training parameters; (square): training set, (circle): testing set. (black): reduced dimension $r=25$, (red): reduced dimension $r=75$.}
\label{fig:param_Pex1}
\end{figure}
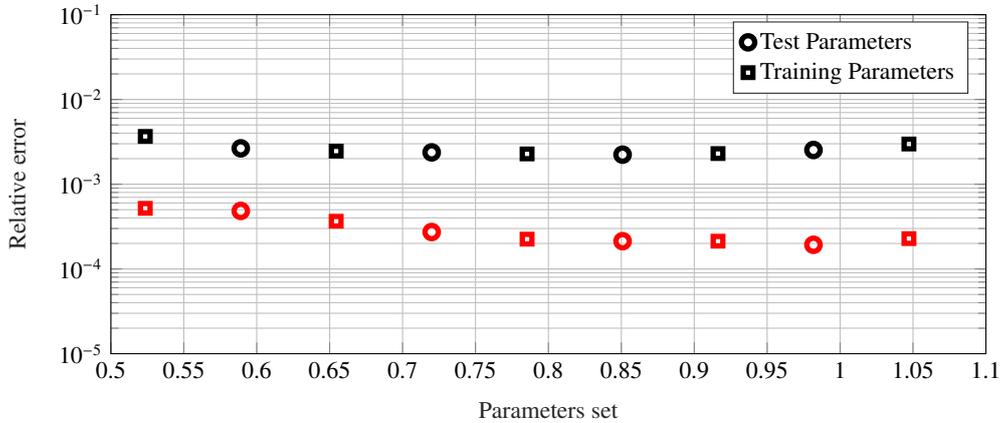

Next, we examine the performance of the parametric non-intrusive ROM on the test set, which consists of the parameters at the midpoint of two successive training parameters. In \figurename{ \ref{fig:param_Pex1}}, we demonstrate the relative errors \eqref{eq:rel_err} of the non-intrusive ROM of orders $ r=25 $ and $ r=75 $ for both test and training parameters. \figurename{ \ref{fig:param_Pex1}} shows that the accuracy of the non-intrusive ROM increases when the reduced dimension increases for training parameters as well as for test parameters.

\begin{figure}[H]
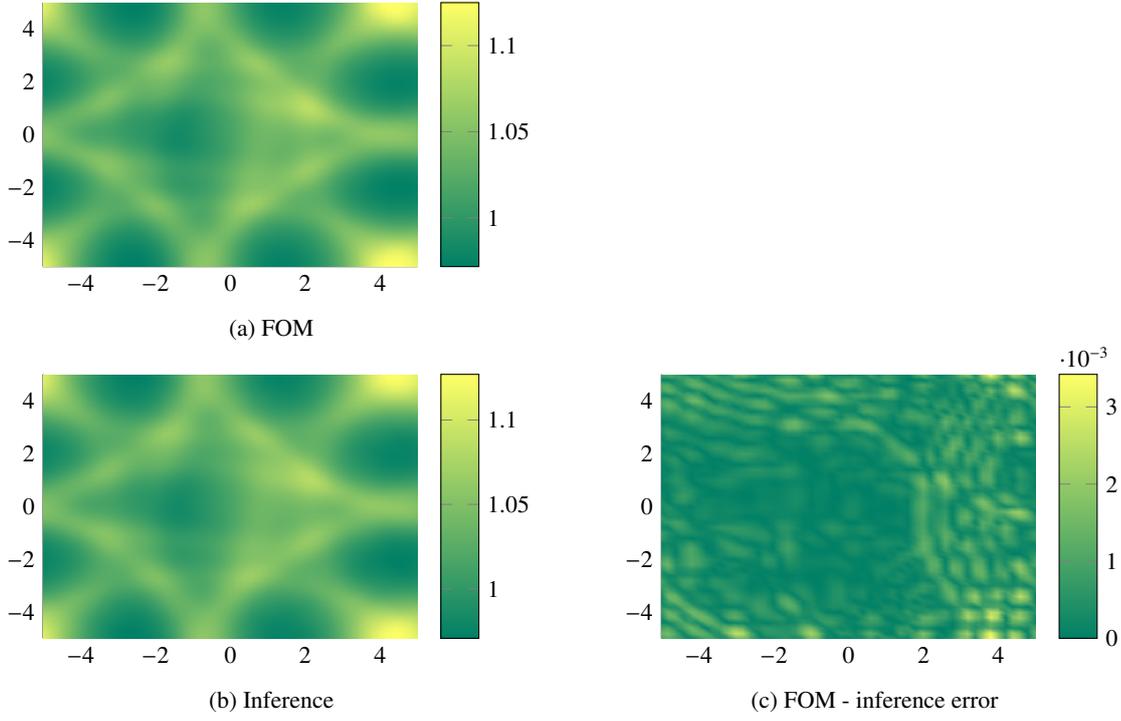

	\begin{subfigure}[t]{0.45\textwidth}
		\centering
		\setlength\fheight{1.5cm}
		\setlength\fwidth{.8\textwidth}
		{\small
    \tikzsetnextfilename{fig//Pex1/fom}%
    \input{fig//Pex1/fom.tikz}%

		}
		\caption{FOM}
		\label{fig:fomPex1}
	\end{subfigure}
\\
	\begin{subfigure}[t]{0.45\textwidth}
		\centering
		\setlength\fheight{1.5cm}
		\setlength\fwidth{.8\textwidth}
		{\small
    \tikzsetnextfilename{fig//Pex1/inf}%
    \input{fig//Pex1/inf.tikz}%

		}
		\caption{Inference}
		\label{fig:InfPex1}
	\end{subfigure}
	\begin{subfigure}[t]{0.45\textwidth}
	\centering
	\setlength\fheight{1.5cm}
	\setlength\fwidth{.8\textwidth}
	{\small
    \tikzsetnextfilename{fig//Pex1/fom_inf}%
    \input{fig//Pex1/fom_inf.tikz}%

	}
	\caption{FOM - inference error}
	\label{fig:fom_infPex1}
\end{subfigure}
	\caption{Single-layer geostrophic adjustment: Comparison of the height field for the parameter $\mu=\frac{5\pi}{24}$  obtained using the FOM and non-intrusive ROM of order $r=75$ at time $T=10$.}
	\label{fig:comparison_HF_param}
\end{figure}

Finally, we show the height field at time $ T=10 $ for the parameter $\mu =  \frac{5\pi}{24}$ and the corresponding absolute errors in \figurename{ \ref{fig:comparison_HF_param}} for the non-intrusive ROM of order $r=75$ which shows that the non-intrusive ROM captures the dynamics of the FOM very well.

\FloatBarrier
\subsection{Single-layer shear instability}
The initial conditions for the second test example on the periodic domain $ [0,10]\times [0,10] $ are given as:
\begin{align*}
&h(x, y, 0) = 1+\Delta h\sin\bigg\{\frac{2\pi}{L}\left[y-\Delta y \sin\left(\frac{2\pi x}{L}\right)\right]\bigg\},\\
&u(x, y, 0) = -\frac{2\pi\Delta h}{\Omega^z L}\cos\bigg\{\frac{2\pi}{L}\left[y-\Delta y \sin\left(\frac{2\pi x}{L}\right)\right]\bigg\},\\
&v(x, y, 0) = -\frac{4\pi^2\Delta h \Delta y}{\Omega^z L^2}\cos\bigg\{\frac{2\pi}{L}\left[y-\Delta y \sin\left(\frac{2\pi x}{L}\right)\right]\bigg\}\cos\left(\frac{2\pi x}{L}\right),
\end{align*}
where $\Delta h=0.2 $, $ \Delta y=0.5 $ and the dimensionless spatial domain length $ L=10 $. This test example illustrates the roll-up of an unstable shear layer \cite{Stewart16}.

\FloatBarrier
\subsubsection{Non-parametric case}

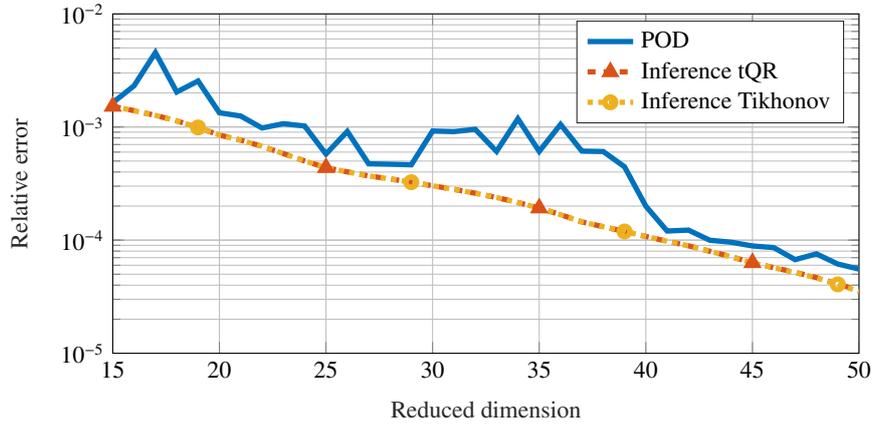
\begin{figure}[H]
		\centering
		\setlength\fheight{1.5cm}
		\setlength\fwidth{.6\textwidth}
		{\small
    \tikzsetnextfilename{fig//Nex2/rel_err}%
    \input{fig//Nex2/rel_err.tikz}%

		}
	\caption{Single-layer shear instability: Relative ROMs errors.}
	\label{fig:rel_errNex2}
\end{figure}

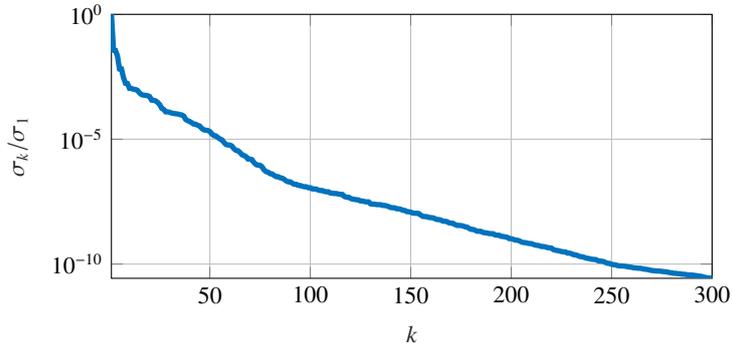
\begin{figure}[H]
	\centering
		\centering
		\setlength\fheight{1.5cm}
		\setlength\fwidth{.5\textwidth}
		{\small
    \tikzsetnextfilename{fig//Nex2/sing}%
    \input{fig//Nex2/sing.tikz}%

		}
		\caption{Single-layer shear instability: Normalized singular values}
		\label{fig:singNex2}
\end{figure}

The performance of the OpInf method is shown in terms of learning the vorticity dynamics for the parameter $\mu= \frac{\pi}{4} $. We collect the snapshots as the discrete state vectors concatenated into $\bw=[\bu,\bv,\bh]^T\in\mathbb{R}^{30000}$. We perform simulations of the FOM \eqref{discrete_fom} in the time domain $[0,60]$, which yields the training set of the size $ \bS_{\text{FOM}}\in\mathbb{R}^{30000\times601} $. The quality of the ROMs in terms of relative errors \eqref{eq:rel_err} are shown in \figurename{ \ref{fig:rel_errNex2}} , which shows that the relative errors of the  non-intrusive ROMs decrease smoothly with decreasing singular values in \figurename{ \ref{fig:singNex2}}, whereas the non-intrusive ROM errors decrease non-smoothly.

The projection error \eqref{eq:proj_err} for $ r=50 $ yields, $\mathcal{E}_{\text{proj}} =3.52\cdot 10^{-5}$. Next, we compare the projection error with relative errors \eqref{eq:rel_err} of ROMs of order $ r=50 $ in \tablename{ \ref{tab:3}}, which again indicates that the non-intrusive ROMs are more accurate.

\begin{table}[H]
	\centering
	\begin{tabular}{|c|c|c|c|}
		\hline
		Method & POD & Non-intrusive (Tikhonov regularizer) & Non-intrusive (tQR)\\\hline
		$\mathcal{E}$ &$5.54\cdot 10^{-5}$ & $ 3.52 \cdot 10^{-5} $  &  $3.58 \cdot 10^{-5} $ \\
		\hline
	\end{tabular}
	\caption{Single-layer shear instability: Comparison of intrusive and non-intrusive ROMs.}
	\label{tab:3}
\end{table}

\begin{figure}[H]
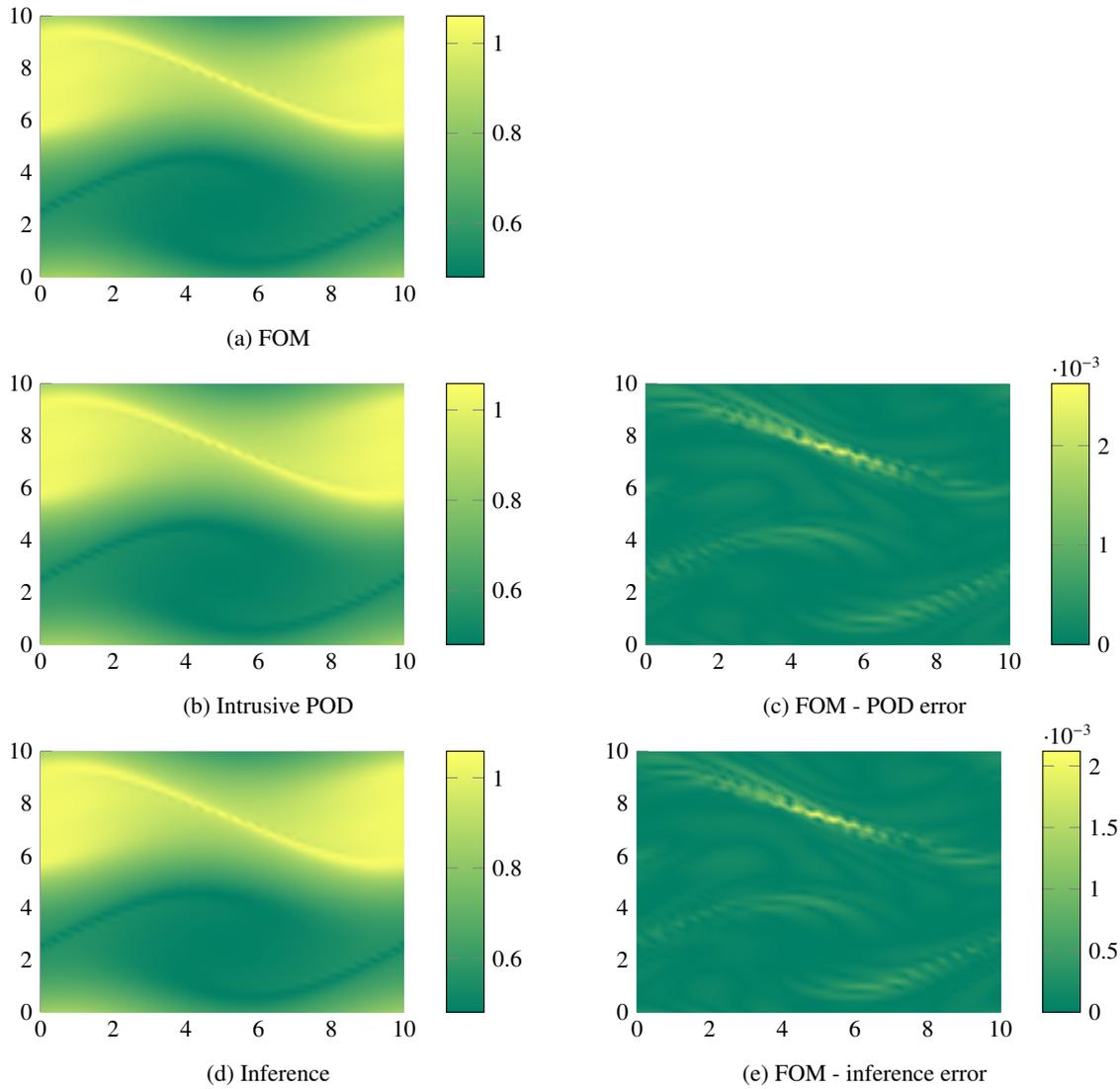

	\begin{subfigure}[t]{0.45\textwidth}
		\centering
		\setlength\fheight{1.5cm}
		\setlength\fwidth{.8\textwidth}
		{\small
    \tikzsetnextfilename{fig//Nex2/fom}%
    \input{fig//Nex2/fom.tikz}%

		}
		\caption{FOM}
		\label{fig:fomNex2}
	\end{subfigure}
	\\
	\begin{subfigure}[t]{0.45\textwidth}
		\centering
		\setlength\fheight{1.5cm}
		\setlength\fwidth{.8\textwidth}
		{\small
    \tikzsetnextfilename{fig//Nex2/pod}%
    \input{fig//Nex2/pod.tikz}%

		}
		\caption{Intrusive POD}
		\label{fig:podNex2}
	\end{subfigure}
	\begin{subfigure}[t]{0.45\textwidth}
		\centering
		\setlength\fheight{1.5cm}
		\setlength\fwidth{.8\textwidth}
		{\small
    \tikzsetnextfilename{fig//Nex2/fom_pod}%
    \input{fig//Nex2/fom_pod.tikz}%

		}
		\caption{FOM -  POD  error}
		\label{fig:fom_podNex2}
	\end{subfigure}
	\\
	\begin{subfigure}[t]{0.45\textwidth}
		\centering
		\setlength\fheight{1.5cm}
		\setlength\fwidth{.8\textwidth}
		{\small
    \tikzsetnextfilename{fig//Nex2/inf}%
    \input{fig//Nex2/inf.tikz}%

		}
		\caption{Inference}
		\label{fig:infNex2}
	\end{subfigure}
	\begin{subfigure}[t]{0.45\textwidth}
		\centering
		\setlength\fheight{1.5cm}
		\setlength\fwidth{.8\textwidth}
		{\small
    \tikzsetnextfilename{fig//Nex2/fom_inf}%
    \input{fig//Nex2/fom_inf.tikz}%

		}
		\caption{FOM - inference error}
		\label{fig:fom_infNex2}
	\end{subfigure}
	\caption{Single-layer shear instability: Comparison of the potential vorticity field obtained using the FOM and ROMs of order $r = 50$ at time $ T=60 $.}
	\label{fig:comparison_VF_t60}
\end{figure}

The potential vorticities of the FOM and ROMs of order $r=50$ as well as corresponding absolute error are shown in \figurename{ \ref{fig:comparison_VF_t60}}, where both the FOM and ROMs share similar roll-up behavior in the vorticity dynamics. In Figures \ref{fig:fom_podNex2}, \ref{fig:fom_infNex2}, we observe that the accuracy of intrusive and non-intrusive ROMs of order $ r=50$ is similar for the vorticity dynamics.

\begin{figure}[H]
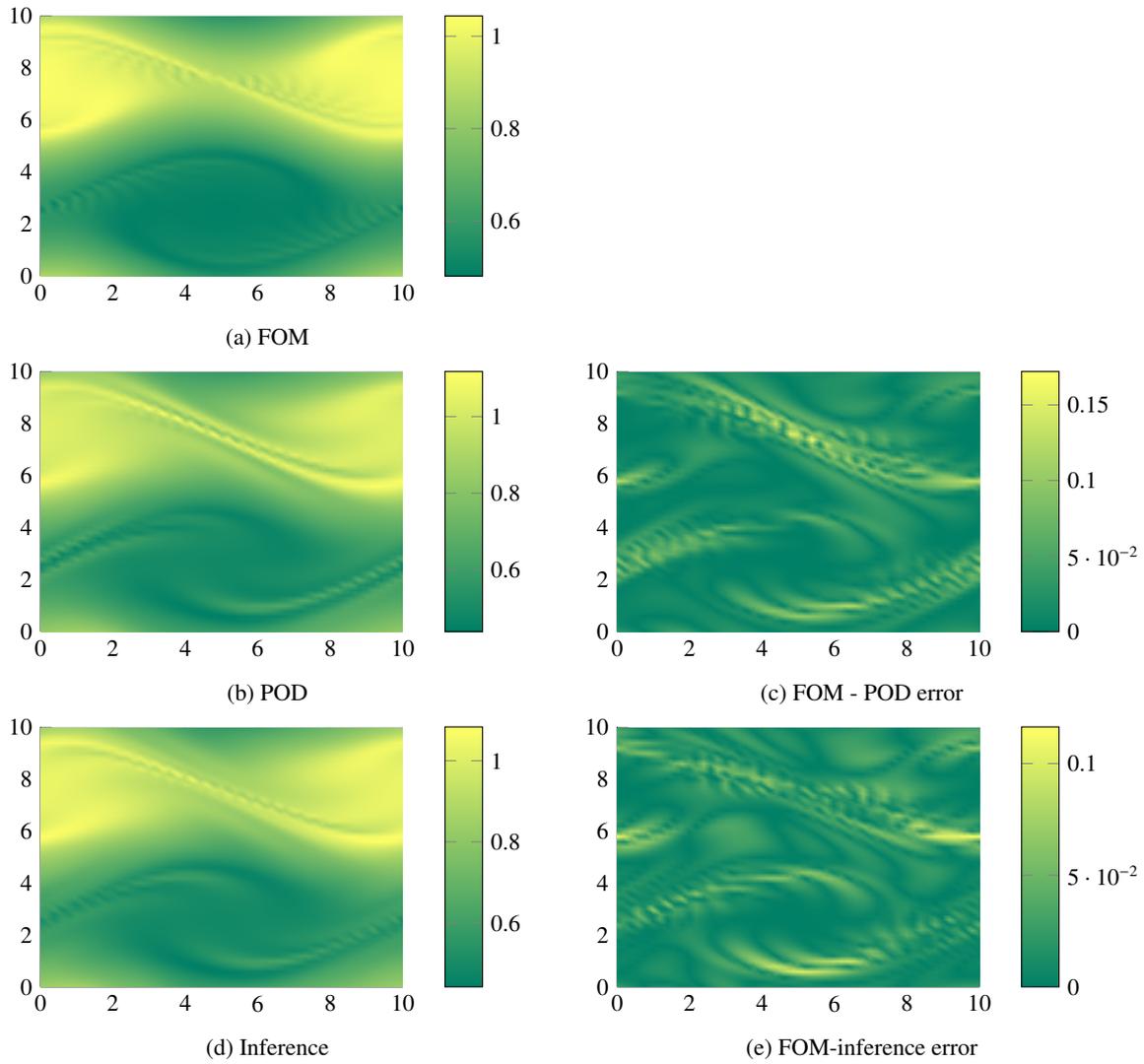

	\begin{subfigure}[t]{0.45\textwidth}
		\centering
		\setlength\fheight{1.5cm}
		\setlength\fwidth{.8\textwidth}
		{\small
    \tikzsetnextfilename{fig//Nex2/Pred/fom}%
    \input{fig//Nex2/Pred/fom.tikz}%

		}
		\caption{FOM}
		\label{fig:fomNex2pred}
	\end{subfigure}
	\\
	\begin{subfigure}[t]{0.45\textwidth}
		\centering
		\setlength\fheight{1.5cm}
		\setlength\fwidth{.8\textwidth}
		{\small
    \tikzsetnextfilename{fig//Nex2/Pred/pod}%
    \input{fig//Nex2/Pred/pod.tikz}%

		}
		\caption{POD}
		\label{fig:podNex2pred}
	\end{subfigure}
	\begin{subfigure}[t]{0.45\textwidth}
		\centering
		\setlength\fheight{1.5cm}
		\setlength\fwidth{.8\textwidth}
		{\small
    \tikzsetnextfilename{fig//Nex2/Pred/fom_pod}%
    \input{fig//Nex2/Pred/fom_pod.tikz}%

		}\\
		\caption{FOM - POD  error}
		\label{fig:fom_podNex2pred}
	\end{subfigure}
	\\
	\begin{subfigure}[t]{0.45\textwidth}
		\centering
		\setlength\fheight{1.5cm}
		\setlength\fwidth{.8\textwidth}
		{\small
    \tikzsetnextfilename{fig//Nex2/Pred/inf}%
    \input{fig//Nex2/Pred/inf.tikz}%

		}
		\caption{Inference}
		\label{fig:infNex2pred}
	\end{subfigure}
	\begin{subfigure}[t]{0.45\textwidth}
		\centering
		\setlength\fheight{1.5cm}
		\setlength\fwidth{.8\textwidth}
		{\small
    \tikzsetnextfilename{fig//Nex2/Pred/fom_inf}%
    \input{fig//Nex2/Pred/fom_inf.tikz}%

		}\\
		\caption{FOM-inference error}
		\label{fig:fom_infNex2pred}
	\end{subfigure}
	\caption{Single-layer shear instability: Prediction of the potential vorticity field obtained using the FOM and ROMs of order $r=50$ at time $T=80$.}
	\label{fig:prediction_VF_t80}
\end{figure}

In \figurename{ \ref{fig:prediction_VF_t80}}, we demonstrate the prediction capability of the ROMs of order $r=50 $ obtained via intrusive POD and non-intrusive OpInf methods by training them in the time interval $ [0,60] $. We set the final time for the prediction as $T=80$. \figurename{ \ref{fig:prediction_VF_t80}} shows that the non-intrusive OpInf solutions are more accurate than the intrusive POD solutions.
% which  indicates that prediction can be performed without accessing the full order operators.

\FloatBarrier
\subsubsection{Parametric Case}
In the last example, we examine the performance of the ROMs for the vorticity dynamics by setting the parameter domain of NTSWE as $ \mathcal{D} = \left[\frac{\pi}{6},\frac{\pi}{3}\right] \subset \mathbb{R} $. In this case, the initial condition is dependent on the angular velocity vector $\Omega^{z}$. The trajectories in the training set are constructed simulating NTSWE \eqref{discrete_fom} on the time domain $[0,T]$  with final time $ T=30 $ and $M=5$ equidistantly distributed parameters  $\mu_1,\mu_2,\ldots,\mu_5 \in \mathcal{D}$. The total size of the concatenated snapshot matrix is $\bS_{\text{FOM}}\in\mathbb{R}^{30000\times 1505} $.

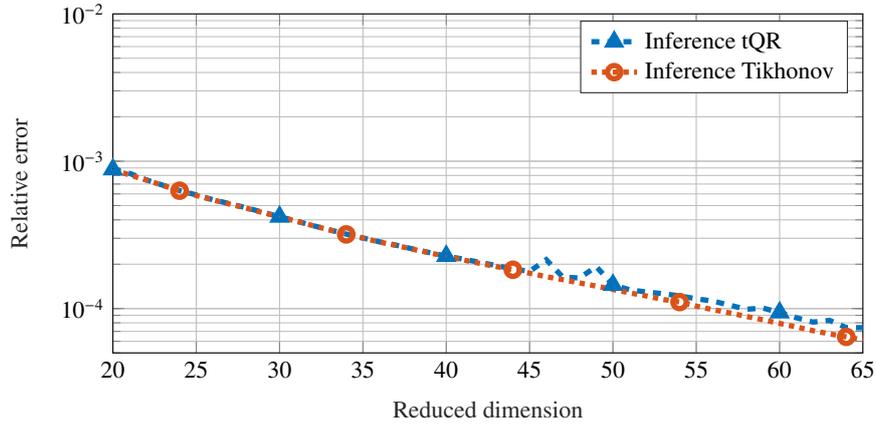
\begin{figure}[H]
		\centering
		\setlength\fheight{4.5cm}
		\setlength\fwidth{.6\textwidth}
		{\small
    \tikzsetnextfilename{fig//Pex2/rel_err_training}%
    \input{fig//Pex2/rel_err_training.tikz}%

		}
	\caption{Single-layer shear instability: Relative ROMs errors.}
	\label{fig:rel_err_Pex2}
\end{figure}

\begin{figure}[H]
	\centering
		\setlength\fheight{1.5cm}
		\setlength\fwidth{.5\textwidth}
		{\small
    \tikzsetnextfilename{fig//Pex2/sing}%
    \input{fig//Pex2/sing.tikz}%

		}
		\caption{Single-layer shear instability: Normalized singular values}
         \label{fig:singPex2}
\end{figure}
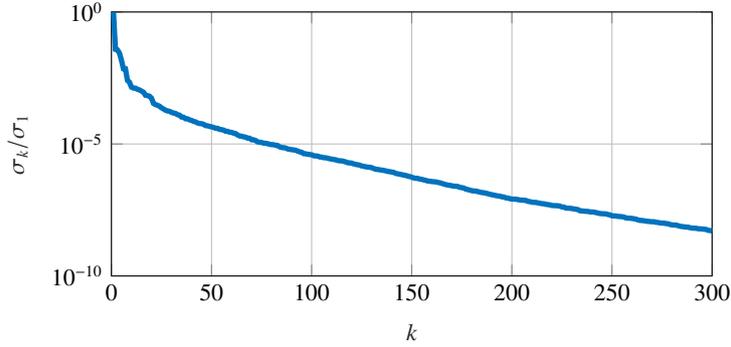

The relative errors \eqref{eq:rel_err} of the non-intrusive ROMs of orders $20$ to $65$ over the training set shown in \figurename{ \ref{fig:rel_err_Pex2}},
are smoothly decreasing when the normalized singular values decrease as shown in \figurename{ \ref{fig:singPex2}} as in the first test example. The projection error \eqref{eq:proj_err} of the FOM for $ r=65 $  is $ \mathcal{E}=6.12\cdot 10^{-5} $. Next, we compare the projection error with the relative error \eqref{eq:rel_err} of the non-intrusive ROMs of order $ r=65 $  in \tablename{ \ref{tab:4}} , which shows that  both regularizers provide accurate simulations of the same order.

\begin{table}[H]
	\centering
	\begin{tabular}{|c|c|c|}
		\hline
		Method  & Non-intrusive (Tikhonov regularizer) & Non-intrusive (tQR)\\\hline
		$\mathcal{E}$  & $ 6.12 \cdot 10^{-5} $  &  $7.41 \cdot 10^{-5} $ \\
		\hline
	\end{tabular}
	\caption{Single-layer shear instability (parametric case): Comparison of non-intrusive ROMs.}
	\label{tab:4}
\end{table}

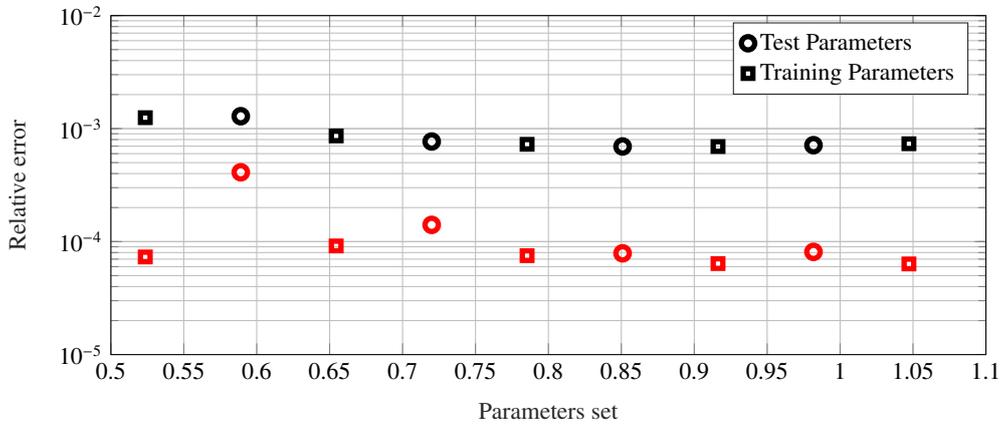
\begin{figure}[H]
	\centering
	\setlength\fheight{4.5cm}
	\setlength\fwidth{.7\textwidth}
	{\small
    \tikzsetnextfilename{fig//Pex2/param_err}%
    \input{fig//Pex2/param_err.tikz}%

	}
	\caption{Single-layer shear instability: Relative error for testing and training parameters; (square):training set, (circle):testing set. (black): reduced dimension $20$, (red): reduced dimension $65$.}
	\label{fig:param_errPex2}
\end{figure}

We again consider the test set consisting of the parameters at the middle points of two successive training parameters. \figurename{ \ref{fig:param_errPex2}} shows the relative errors for parameters in the test and training sets for non-intrusive ROMs of order $20$ and $65$. In \figurename{ \ref{fig:param_errPex2}}, the inferred solution for the test parameter $ \mu=\frac{9\pi}{48} $ is less accurate than at other parameters in the test set for reduced dimension $65$. This indicates that the truncation tolerance for the QR-CP method can deteriorate the accuracy of the ROMs.

\begin{figure}[H]
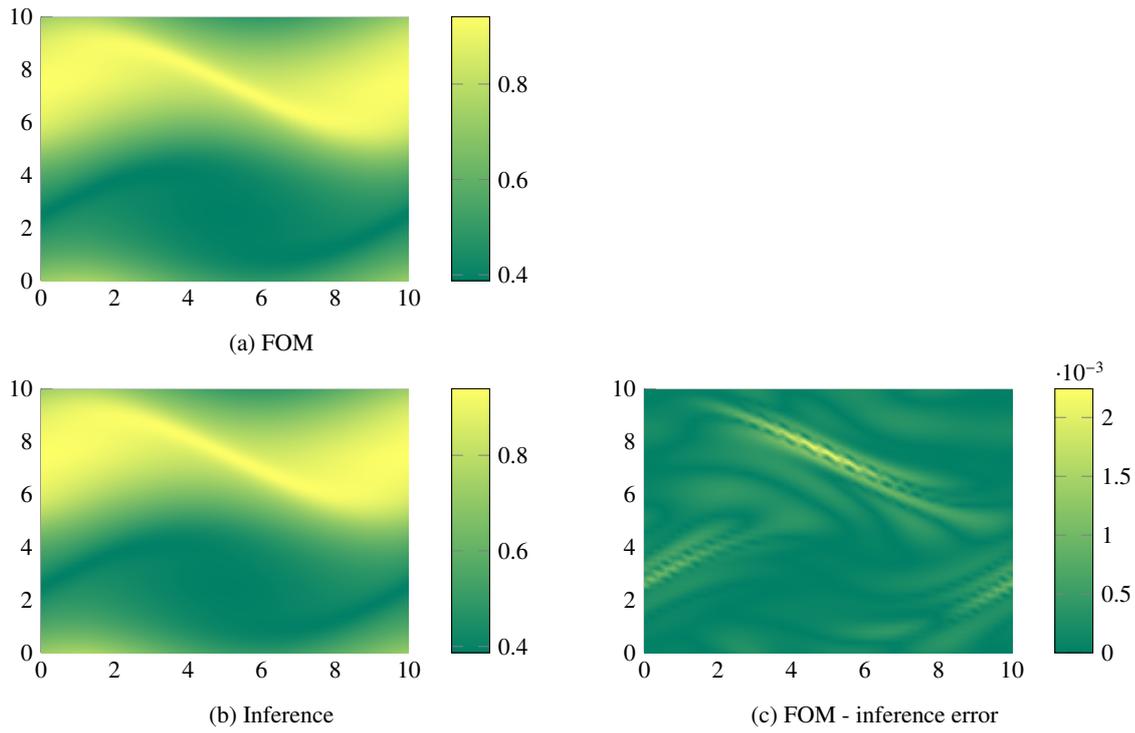

	\begin{subfigure}[t]{0.45\textwidth}
		\centering
		\setlength\fheight{1.5cm}
		\setlength\fwidth{.8\textwidth}
		{\small
    \tikzsetnextfilename{fig//Pex2/fom}%
    \input{fig//Pex2/fom.tikz}%

		}
		\caption{FOM}
		\label{fig:fomPex2}
	\end{subfigure}
	\\
	\begin{subfigure}[t]{0.45\textwidth}
		\centering
		\setlength\fheight{1.5cm}
		\setlength\fwidth{.8\textwidth}
		{\small
    \tikzsetnextfilename{fig//Pex2/inf}%
    \input{fig//Pex2/inf.tikz}%

		}
		\caption{Inference}
		\label{fig:infPex2}
	\end{subfigure}
	\begin{subfigure}[t]{0.45\textwidth}
		\centering
		\setlength\fheight{1.5cm}
		\setlength\fwidth{.8\textwidth}
		{\small
    \tikzsetnextfilename{fig//Pex2/fom_inf}%
    \input{fig//Pex2/fom_inf.tikz}%

		}
		\caption{FOM - inference error}
		\label{fig:fom_infPex2}
	\end{subfigure}
	\caption{Single-layer shear instability: Comparison of the potential vorticity field for the parameter $\mu=\frac{5\pi}{24}$  obtained using the FOM and non-intrusive ROM of order $r=65$ at time $T=30$.}
	\label{fig:comparison_VF_param}
\end{figure}

Finally, we show the potential vorticity of the FOM, the non-intrusive ROM of order $r=65$ and the corresponding absolute error at time $ T=30 $ for the parameter $\mu=\frac{5\pi}{24}$ in \figurename{ \ref{fig:comparison_VF_param}}, which shows that the non-intrusive ROM captures the dynamics of the NTSWE accurately.

\clearpage
\section{Conclusions}
We have constructed  data-driven projection-based ROMs of NTSWE by exploiting the structure of the equations. The OpInf framework is used to construct the non-intrusive ROM. Since the least-squares problem of the OpInf method may suffer from  ill-conditioning, the solutions are regularized using the QR factorization as an alternative to Tikhonov regularization. The performance of the inferred models is examined in terms of prediction capabilities on two test examples. Numerical results show that the prediction of the learned model of OpInf is more accurate than the intrusive POD.

\section*{Acknowledgments}
This work was supported by 100/2000 Ph.D. Scholarship Program of the Turkish Higher Education Council. The first author thanks for the hospitality of the Max Planck Institute   for Dynamics of Complex Technical Systems, Magdeburg, Germany.

%%%%%%%%%%%%%%%%%%%%%%%%%%%%%%%%%%%%%%%%%%%%%%%%%%%%%%%%%%%%%%%%%%%%%%%%%%%%%%%%%%%%%%%%%%%%
%\bibliographystyle{plain}
%\bibliography{references}

\end{document}

%% file: fig/Nex1/sing.tikz
\definecolor{mycolor1}{rgb}{0.00000,0.44700,0.74100}%
\begin{tikzpicture}

\begin{axis}[%
width=0.996\fwidth,
height=\fheight,
at={(0\fwidth,0\fheight)},
scale only axis,
xmin=1,
xmax=300,
xlabel style={font=\color{white!15!black}},
xlabel={$k$},
ymode=log,
ymin=1.48673616355391e-07,
ymax=1,
yminorticks=true,
ylabel style={font=\color{white!15!black}},
ylabel={$\sigma_k / \sigma_1$},
axis background/.style={fill=white},
xmajorgrids,
xminorgrids,
ymajorgrids,
yminorgrids
]
\addplot [color=mycolor1, line width=2.0pt, forget plot]
  table[row sep=crcr]{%
1	1\\
2	0.0367422492261725\\
3	0.0358705625819637\\
4	0.0205161371801504\\
5	0.0202558635445242\\
6	0.0179676651142923\\
7	0.0177660789201208\\
8	0.0162204977066392\\
9	0.0161103060197269\\
10	0.0139730440778083\\
11	0.0138930431731339\\
12	0.0123295248575106\\
13	0.0118563668497187\\
14	0.0104908011659998\\
15	0.0104575929345425\\
16	0.00831973385057911\\
17	0.0082713480393732\\
18	0.00801758385430299\\
19	0.00778140698221198\\
20	0.00622727443341818\\
21	0.00621436255939697\\
22	0.00588724514231612\\
23	0.00581585511043447\\
24	0.00557923387765326\\
25	0.00523288088472439\\
26	0.0047361104121916\\
27	0.00460275893786423\\
28	0.00319747328915245\\
29	0.00314052322433195\\
30	0.00275906108339724\\
31	0.00270173581592119\\
32	0.00214726907334845\\
33	0.00210527954944525\\
34	0.00204414777942131\\
35	0.00198432181649521\\
36	0.00182967487666267\\
37	0.00180133880249858\\
38	0.00175715981157368\\
39	0.00175126798332458\\
40	0.00153840976329317\\
41	0.00152933911134942\\
42	0.00144832830914156\\
43	0.00144374600659114\\
44	0.0014086767990138\\
45	0.00140554749786361\\
46	0.00130955694222922\\
47	0.00130699683878788\\
48	0.00127135700498759\\
49	0.00126495117852291\\
50	0.00122311507028965\\
51	0.0012073496589327\\
52	0.00114032516444624\\
53	0.00111055282962288\\
54	0.00108695176553373\\
55	0.00107616848842872\\
56	0.00103960244404627\\
57	0.00100263112159857\\
58	0.00094464942932017\\
59	0.000933934228169451\\
60	0.000898808593249485\\
61	0.000893821380512366\\
62	0.000826784949339561\\
63	0.000803749997352407\\
64	0.000794663169659829\\
65	0.000771206629099252\\
66	0.000747479086298756\\
67	0.00073789108105363\\
68	0.000650388111366276\\
69	0.000639257197784499\\
70	0.000626719717906782\\
71	0.000605088202393837\\
72	0.000549646739569798\\
73	0.000547829025293216\\
74	0.000509345716951684\\
75	0.00048865800785222\\
76	0.000481351349599745\\
77	0.000471149682732272\\
78	0.000442554485569994\\
79	0.000439634913273518\\
80	0.000435572798593394\\
81	0.000406639269332033\\
82	0.000402971695950094\\
83	0.000401924325770426\\
84	0.000398112037288884\\
85	0.000394000849929547\\
86	0.000372767035345977\\
87	0.000364786064133489\\
88	0.000362129855895094\\
89	0.000348486846512083\\
90	0.000344915348499607\\
91	0.000332308318755824\\
92	0.000325072393031455\\
93	0.000320222066285685\\
94	0.000308401893841433\\
95	0.000299446831371625\\
96	0.0002875561681165\\
97	0.00027546901356471\\
98	0.000274701463853827\\
99	0.000272852052346776\\
100	0.000259771712369707\\
101	0.000255258387724718\\
102	0.000241206767868689\\
103	0.000236203698619111\\
104	0.000234592981614891\\
105	0.000228315333686429\\
106	0.000218370545082734\\
107	0.0002071060270269\\
108	0.000181772501900792\\
109	0.000178536818037442\\
110	0.00017109468806681\\
111	0.00016888433936246\\
112	0.000166860522362632\\
113	0.000163808869652574\\
114	0.000157759918962346\\
115	0.000156915476030351\\
116	0.000145903165528143\\
117	0.000140838347278993\\
118	0.000135350665289309\\
119	0.000134903863082515\\
120	0.000128014562721614\\
121	0.000126570098298895\\
122	0.000124455391766348\\
123	0.000121187224412126\\
124	0.000120548693927\\
125	0.000114408338190018\\
126	0.000108899583734986\\
127	0.00010382198135644\\
128	0.000102013091544116\\
129	9.93168883429022e-05\\
130	9.703065892116e-05\\
131	9.45094228949377e-05\\
132	9.43332282370615e-05\\
133	9.09297484441696e-05\\
134	8.71630782199355e-05\\
135	7.9732594369064e-05\\
136	7.60549145567286e-05\\
137	7.49979457596861e-05\\
138	7.29302324273247e-05\\
139	7.20527963525429e-05\\
140	7.03821702935371e-05\\
141	6.75446559052377e-05\\
142	6.02029944513942e-05\\
143	5.85011267864702e-05\\
144	5.82905773952683e-05\\
145	5.44521571337621e-05\\
146	5.3753312211748e-05\\
147	5.22193928106134e-05\\
148	5.18593972870161e-05\\
149	4.74284526626802e-05\\
150	4.72490446287016e-05\\
151	4.61220520887077e-05\\
152	4.37348202725787e-05\\
153	4.34347941914639e-05\\
154	4.09815184747782e-05\\
155	4.07165471345582e-05\\
156	3.78173748646182e-05\\
157	3.64379557535029e-05\\
158	3.5052755448755e-05\\
159	3.46235640484925e-05\\
160	3.33348379973858e-05\\
161	3.24412693401617e-05\\
162	3.17817386082847e-05\\
163	2.82648112539588e-05\\
164	2.80012185350173e-05\\
165	2.75263887406937e-05\\
166	2.66253712216009e-05\\
167	2.52972629429615e-05\\
168	2.38878106304738e-05\\
169	2.35638725160701e-05\\
170	2.25506163467677e-05\\
171	2.14507827448185e-05\\
172	1.9627687896208e-05\\
173	1.91603342549336e-05\\
174	1.82871682890939e-05\\
175	1.72124610228607e-05\\
176	1.67036830977265e-05\\
177	1.59190405144738e-05\\
178	1.58224792858126e-05\\
179	1.46224550940649e-05\\
180	1.44797534449031e-05\\
181	1.36917541476168e-05\\
182	1.36297023415463e-05\\
183	1.35183872209186e-05\\
184	1.28653586815669e-05\\
185	1.23810172302864e-05\\
186	1.14389043606143e-05\\
187	1.13119717369181e-05\\
188	1.11315154808682e-05\\
189	1.10067874137595e-05\\
190	1.01322364438189e-05\\
191	9.80938519812449e-06\\
192	9.68524951595556e-06\\
193	9.21339999816814e-06\\
194	9.16701756260226e-06\\
195	8.22934677368045e-06\\
196	8.09085458302382e-06\\
197	8.01399001942069e-06\\
198	7.82472142751122e-06\\
199	7.66966574652198e-06\\
200	7.02154239594281e-06\\
201	6.74545552045937e-06\\
202	6.44815855581498e-06\\
203	6.37178425738906e-06\\
204	6.19841925903731e-06\\
205	5.9327133464681e-06\\
206	5.37016749459239e-06\\
207	5.10114501957994e-06\\
208	5.07079299718335e-06\\
209	4.76060096090861e-06\\
210	4.7178300910278e-06\\
211	4.44167360941019e-06\\
212	4.26814120339398e-06\\
213	3.96748822994984e-06\\
214	3.84229026842733e-06\\
215	3.81972656789631e-06\\
216	3.72872352254919e-06\\
217	3.68705046973865e-06\\
218	3.44548083239703e-06\\
219	3.37402294428573e-06\\
220	3.12928960906642e-06\\
221	3.0286259362544e-06\\
222	2.97971379556408e-06\\
223	2.76986831875772e-06\\
224	2.74455980263278e-06\\
225	2.6576511429425e-06\\
226	2.58906222445561e-06\\
227	2.52470269920861e-06\\
228	2.46956487222875e-06\\
229	2.28575029971868e-06\\
230	2.14163484609224e-06\\
231	2.07751571150066e-06\\
232	2.05471568797017e-06\\
233	2.02501940264663e-06\\
234	1.83257993684473e-06\\
235	1.79136711417591e-06\\
236	1.70659002570133e-06\\
237	1.69431675661871e-06\\
238	1.64075624503245e-06\\
239	1.49122829018971e-06\\
240	1.4817546699574e-06\\
241	1.38669844400121e-06\\
242	1.34720548019206e-06\\
243	1.26208912919927e-06\\
244	1.24860759876268e-06\\
245	1.21002938335218e-06\\
246	1.18218581791885e-06\\
247	1.16733046855811e-06\\
248	1.08003807401793e-06\\
249	1.06223161242256e-06\\
250	1.02887678584553e-06\\
251	1.02023583969957e-06\\
252	9.49260626108198e-07\\
253	9.0233492824464e-07\\
254	8.57716989586391e-07\\
255	8.21433347088091e-07\\
256	8.09840238803659e-07\\
257	7.55268372978531e-07\\
258	7.44040478957544e-07\\
259	7.01279342911431e-07\\
260	6.97092296202505e-07\\
261	6.68205784190853e-07\\
262	6.54729632812146e-07\\
263	6.47905641227234e-07\\
264	5.86644623890192e-07\\
265	5.71711268644686e-07\\
266	5.46829139995729e-07\\
267	5.40165168167203e-07\\
268	5.01648725602164e-07\\
269	4.77342063516095e-07\\
270	4.65145592418995e-07\\
271	4.52677761600441e-07\\
272	4.43688227550446e-07\\
273	4.09761044488256e-07\\
274	4.03967829152852e-07\\
275	3.78811255940286e-07\\
276	3.61222403285384e-07\\
277	3.51215929304838e-07\\
278	3.47993114421711e-07\\
279	3.42223298444749e-07\\
280	3.2116999204199e-07\\
281	3.1195662903824e-07\\
282	2.90376009444714e-07\\
283	2.82692576463673e-07\\
284	2.80336798624672e-07\\
285	2.77142086110404e-07\\
286	2.57350390013313e-07\\
287	2.34836779009933e-07\\
288	2.33908461546776e-07\\
289	2.25173453408269e-07\\
290	2.20872011094234e-07\\
291	2.07847538476863e-07\\
292	2.01655041514042e-07\\
293	1.9037145714497e-07\\
294	1.87951159608933e-07\\
295	1.86086733487961e-07\\
296	1.74093795675939e-07\\
297	1.6918007135697e-07\\
298	1.55113794813545e-07\\
299	1.51116790063272e-07\\
300	1.48673616355391e-07\\
};
\end{axis}
\end{tikzpicture}%

%% file: fig/tikL.tikz
\begin{tikzpicture}

\begin{axis}[%
width=\fwidth,
height=0.944\fheight,
at={(0\fwidth,0\fheight)},
scale only axis,
xmode=log,
xmin=1.21662233407247e-06,
xmax=0.00832049655086458,
xminorticks=true,
xlabel style={font=\color{white!15!black}},
xlabel={$\|\mathcal{A}\mathcal{X}-\dot{\hat{\mathbf{S}}}^T\|_{2}^{2}$},
ymode=log,
ymin=0.01,
ymax=9.24152382922357,
yminorticks=true,
ylabel style={font=\color{white!15!black}},
ylabel={$\|\mathcal{X}\|_{2}^{2}$},
axis background/.style={fill=white},
xmajorgrids,
xminorgrids,
ymajorgrids,
yminorgrids
]
\addplot [color=red, line width=2.0pt, draw=none, mark=o, mark options={solid, red}, forget plot]
  table[row sep=crcr]{%
0.00832049655086458	0.0181559215105904\\
0.00356111195115413	0.0184765537467213\\
0.00174349541658423	0.0212342267919317\\
0.0005826489512117	0.0577713629513265\\
0.000128898983989283	0.189160901013724\\
3.04166764431655e-05	0.457872551657628\\
9.58769349229179e-06	1.09059598406734\\
3.01339218008378e-06	3.30041846675714\\
1.21662233407247e-06	9.24152382922357\\
};
\addplot [color=blue, line width=2.0pt, draw=none, mark=asterisk, mark options={solid, blue}, forget plot]
  table[row sep=crcr]{%
0.0005826489512117	0.0577713629513265\\
};
\end{axis}
\end{tikzpicture}%

%% file: fig/tqrL.tikz
\begin{tikzpicture}

\begin{axis}[%
width=\fwidth,
height=0.922\fheight,
at={(0\fwidth,0\fheight)},
scale only axis,
xmode=log,
xmin=1.54273680757627e-07,
xmax=2.02142589166756,
xminorticks=true,
xlabel style={font=\color{white!15!black}},
xlabel={$\|\mathcal{A}\mathcal{X}-\dot{\hat{\mathbf{S}}}^T\|_{2}^{2}$},
ymode=log,
ymin=0.001,
ymax=399.111468418119,
yminorticks=true,
ylabel style={font=\color{white!15!black}},
ylabel={$\|\mathcal{X}\|_{2}^{2}$},
axis background/.style={fill=white},
xmajorgrids,
xminorgrids,
ymajorgrids,
yminorgrids
]
\addplot [color=red, line width=2.0pt, draw=none, mark=o, mark options={solid, red}, forget plot]
  table[row sep=crcr]{%
2.02142589166756	0.017814320525726\\
0.0117150826644706	0.0181738267326114\\
0.00357678681781457	0.0200418584863245\\
0.000279419174280692	0.223479000712692\\
1.47584283040925e-05	1.63130635286631\\
1.29617434435622e-06	16.5741323086134\\
1.54273680757627e-07	399.111468418119\\
};
\addplot [color=blue, line width=2.0pt, draw=none, mark=asterisk, mark options={solid, blue}, forget plot]
  table[row sep=crcr]{%
0.00357678681781457	0.0200418584863245\\
};
\end{axis}
\end{tikzpicture}%

%% file: fig/Nex1/rel_err.tikz
\definecolor{mycolor1}{rgb}{0.00000,0.44700,0.74100}%
\definecolor{mycolor2}{rgb}{0.85000,0.32500,0.09800}%
\definecolor{mycolor3}{rgb}{0.92900,0.69400,0.12500}%
\begin{tikzpicture}

\begin{axis}[%
width=\fwidth,
height=0.999\fheight,
at={(0\fwidth,0\fheight)},
scale only axis,
xmin=20,
xmax=75,
xlabel style={font=\color{white!15!black}},
xlabel={Reduced dimension},
ymode=log,
ymin=0.001,
ymax=0.1,
yminorticks=true,
ylabel style={font=\color{white!15!black}},
ylabel={Relative error},
axis background/.style={fill=white},
xmajorgrids,
ymajorgrids,
yminorgrids,
legend style={legend cell align=left, align=left, draw=white!15!black}
]
\addplot [color=mycolor1, line width=2.0pt]
  table[row sep=crcr]{%
20	0.0231222798222139\\
21	0.0223348598733279\\
22	0.0224665990066849\\
23	0.0233653190720527\\
24	0.0224963572274093\\
25	0.0162428785808979\\
26	0.019217722532397\\
27	0.013899426276548\\
28	0.0161201149782007\\
29	0.0130546616561603\\
30	0.0129424362858331\\
31	0.0126365861653282\\
32	0.012470648813741\\
33	0.0161022768697426\\
34	0.0125810186131032\\
35	0.0128731169567085\\
36	0.0127078215224835\\
37	0.0103902570560437\\
38	0.0104779040264948\\
39	0.010248083986997\\
40	0.0104102674162259\\
41	0.0102070295972362\\
42	0.0104656354389395\\
43	0.012611067321742\\
44	0.011929552206782\\
45	0.011136832242774\\
46	0.0111371057884613\\
47	0.00968204723224204\\
48	0.00985543784529435\\
49	0.00932701043866755\\
50	0.00901553294575273\\
51	0.00788314299198697\\
52	0.00864474421405988\\
53	0.00787294422029295\\
54	0.00838872392072236\\
55	0.00708200965051037\\
56	0.00782310992107259\\
57	0.00679356689600884\\
58	0.0067835408693295\\
59	0.00655876539646655\\
60	0.00728979388882288\\
61	0.00768028416230808\\
62	0.00697943821379672\\
63	0.00548398324690065\\
64	0.00615028156832754\\
65	0.00509658196822498\\
66	0.00613824428626037\\
67	0.00465185855176062\\
68	0.00536882067812836\\
69	0.00434375183176517\\
70	0.00511118424461933\\
71	0.00352914873627098\\
72	0.00354254602303051\\
73	0.00378107698138684\\
74	0.00505080778509693\\
75	0.00357074632608655\\
};
\addlegendentry{POD}

\addplot [color=mycolor2, line width=2.0pt,dashed, mark size=2.0pt, mark=triangle, mark repeat=10, mark options={solid, mycolor2}]
  table[row sep=crcr]{%
20	0.017765818874413\\
21	0.0166505777906084\\
22	0.0155818427065828\\
23	0.0144629312738844\\
24	0.0133506058765584\\
25	0.0122891393801748\\
26	0.0113458907478222\\
27	0.0103765759497528\\
28	0.00987481746142905\\
29	0.00936532585275001\\
30	0.00895228460831372\\
31	0.00853748328018912\\
32	0.0082647446250331\\
33	0.00799380342299093\\
34	0.00772967633442402\\
35	0.00747224241746755\\
36	0.0072461819725278\\
37	0.00702012532922626\\
38	0.00679804495532205\\
39	0.00657002447885008\\
40	0.00638850358374686\\
41	0.00620390132043491\\
42	0.00603353583611242\\
43	0.0058593421525912\\
44	0.00568855337650408\\
45	0.00551327132602745\\
46	0.0053564651576375\\
47	0.0051955650974235\\
48	0.00503859409646814\\
49	0.00487822860891225\\
50	0.00472337057697306\\
51	0.00456743925016345\\
52	0.00442369160762358\\
53	0.00428290484155747\\
54	0.00414355458380781\\
55	0.00400224409051562\\
56	0.0038657202393461\\
57	0.00373425512359576\\
58	0.00361355145026339\\
59	0.00349153999289376\\
60	0.00337460474590883\\
61	0.00325483001420364\\
62	0.00314873201237591\\
63	0.00304506933015278\\
64	0.00294020749011128\\
65	0.00283790428133426\\
66	0.00273834719950099\\
67	0.00263767648173775\\
68	0.00255672816443437\\
69	0.00247602089144846\\
70	0.00239589810970148\\
71	0.00231870601982473\\
72	0.00225304061136867\\
73	0.00218585854471282\\
74	0.00212608128038712\\
75	0.00206951016403458\\
};
\addlegendentry{Inference tQR}

\addplot [color=mycolor3, line width=2.0pt, dotted, mark size=2.0pt, mark=o, mark repeat=10, mark phase = 5, mark options={solid, mycolor3}]
  table[row sep=crcr]{%
20	0.0177658157585368\\
21	0.0166505732699497\\
22	0.0155818395795127\\
23	0.0144629258537384\\
24	0.0133506038448911\\
25	0.0122891367117734\\
26	0.0113458882369081\\
27	0.0103765731814828\\
28	0.00987481040295435\\
29	0.00936532380124943\\
30	0.00895228241641019\\
31	0.00853748118937196\\
32	0.00826474410822739\\
33	0.00799381465128664\\
34	0.00772967458129579\\
35	0.00747224165532344\\
36	0.00724618104443582\\
37	0.00702012433748336\\
38	0.00679804398401603\\
39	0.00657002306519471\\
40	0.0063885021730525\\
41	0.00620389988387328\\
42	0.0060335343450628\\
43	0.0058593404592915\\
44	0.00568855216735708\\
45	0.00551326887724357\\
46	0.00535646091247588\\
47	0.00519556262847299\\
48	0.0050385909868616\\
49	0.00487822459206674\\
50	0.00472337216781822\\
51	0.00456743080519498\\
52	0.00442368900846473\\
53	0.00428290009441532\\
54	0.00414354793255351\\
55	0.00400224018846463\\
56	0.00386571563455351\\
57	0.00373425105217319\\
58	0.00361354689732895\\
59	0.00349153516770206\\
60	0.00337459592576019\\
61	0.00325482143256481\\
62	0.00314872453225251\\
63	0.00304506152004164\\
64	0.00294019878266377\\
65	0.00283789417539528\\
66	0.00273830832807712\\
67	0.00263764599163119\\
68	0.00255670811160236\\
69	0.00247600560245491\\
70	0.00239587691504354\\
71	0.00231869181109077\\
72	0.00225301270252253\\
73	0.00218581414039884\\
74	0.00212601402912462\\
75	0.00206944635286579\\
};
\addlegendentry{Inference Tikhonov}

\end{axis}
\end{tikzpicture}%

%% file: fig/Pex1/sing.tikz
\definecolor{mycolor1}{rgb}{0.00000,0.44700,0.74100}%
\begin{tikzpicture}

\begin{axis}[%
width=0.961\fwidth,
height=\fheight,
at={(0\fwidth,0\fheight)},
scale only axis,
xmin=0,
xmax=300,
xlabel style={font=\color{white!15!black}},
xlabel={$k$},
ymode=log,
ymin=1e-10,
ymax=1,
yminorticks=true,
ylabel style={font=\color{white!15!black}},
ylabel={$\sigma_k / \sigma_1$},
axis background/.style={fill=white},
xmajorgrids,
xminorgrids,
ymajorgrids,
yminorgrids
]
\addplot [color=mycolor1, line width=2.0pt, forget plot]
  table[row sep=crcr]{%
1	1\\
2	0.0433725918279657\\
3	0.0411766566499045\\
4	0.0298029771056363\\
5	0.0286625663110297\\
6	0.0201025830192869\\
7	0.0172889333700223\\
8	0.0105977752489008\\
9	0.00984658005647605\\
10	0.00524563229155937\\
11	0.00425096992188861\\
12	0.00415167259297337\\
13	0.00368238684243489\\
14	0.00355014642559343\\
15	0.0029227921529007\\
16	0.00277207655339056\\
17	0.00243893139025475\\
18	0.00223048743059324\\
19	0.00166277986325163\\
20	0.00159808002880038\\
21	0.00155980362482298\\
22	0.00145769567594738\\
23	0.00123961102270861\\
24	0.00113607750753259\\
25	0.00101974203909917\\
26	0.000992905925266757\\
27	0.000897964135490364\\
28	0.000807850748542277\\
29	0.00079968237573173\\
30	0.00077964322900921\\
31	0.000673836938745789\\
32	0.000634507889759658\\
33	0.000611047567472915\\
34	0.000548417252818755\\
35	0.000500051000536135\\
36	0.000488314532302112\\
37	0.000469942343704369\\
38	0.000424112408382182\\
39	0.000421532367907762\\
40	0.000381327650015379\\
41	0.000373143209595761\\
42	0.000365768961296602\\
43	0.000336014342356163\\
44	0.000303185842679825\\
45	0.000292515210084426\\
46	0.00028188828783468\\
47	0.000261444295308191\\
48	0.00025404996717056\\
49	0.000243475112999903\\
50	0.000227830187388024\\
51	0.000215222659951307\\
52	0.00020868624278303\\
53	0.00020250784835371\\
54	0.00018719918507574\\
55	0.000175906486900879\\
56	0.000173286038576881\\
57	0.000159749628669572\\
58	0.000149220398430443\\
59	0.000145897842680726\\
60	0.000132148972381899\\
61	0.000127863556815996\\
62	0.00012182990965247\\
63	0.000115294058627006\\
64	0.000108266251185709\\
65	0.000102260679803315\\
66	9.54612508778814e-05\\
67	9.36941779333682e-05\\
68	8.98905496002607e-05\\
69	8.53171097023718e-05\\
70	8.05578994476474e-05\\
71	7.80823641858363e-05\\
72	7.34800771946526e-05\\
73	7.10642007273689e-05\\
74	6.9122156106783e-05\\
75	6.53260707202963e-05\\
76	6.10030597648275e-05\\
77	5.86466141963761e-05\\
78	5.43281033985863e-05\\
79	5.24603311741585e-05\\
80	5.09612324087232e-05\\
81	4.7165410662256e-05\\
82	4.42474291978575e-05\\
83	4.28943429292145e-05\\
84	4.2727799970312e-05\\
85	4.07770035078603e-05\\
86	3.85920926109211e-05\\
87	3.63696024506789e-05\\
88	3.54467629401617e-05\\
89	3.30159802752848e-05\\
90	3.12732040093065e-05\\
91	3.00927304774105e-05\\
92	2.99349878544901e-05\\
93	2.91223664368736e-05\\
94	2.81077993405482e-05\\
95	2.6073945517583e-05\\
96	2.5088047650999e-05\\
97	2.48659793496023e-05\\
98	2.42810487007363e-05\\
99	2.36681214325388e-05\\
100	2.17228658285932e-05\\
101	2.06497450549369e-05\\
102	2.03360945740911e-05\\
103	1.94381101779359e-05\\
104	1.83650382974746e-05\\
105	1.79321449055195e-05\\
106	1.67581937062446e-05\\
107	1.64578548044234e-05\\
108	1.615535347525e-05\\
109	1.53303129951677e-05\\
110	1.49333979376096e-05\\
111	1.47546104593724e-05\\
112	1.40988311119735e-05\\
113	1.36734613260455e-05\\
114	1.28782304115867e-05\\
115	1.27120554994865e-05\\
116	1.16388419438619e-05\\
117	1.14718095432036e-05\\
118	1.14132421711438e-05\\
119	1.10678832463831e-05\\
120	1.0478120342679e-05\\
121	1.01985629162439e-05\\
122	9.69686194909429e-06\\
123	9.26896179692478e-06\\
124	9.15676111051008e-06\\
125	8.63645650822239e-06\\
126	8.4773051382101e-06\\
127	7.73503759673962e-06\\
128	7.60024291397106e-06\\
129	7.39719995036009e-06\\
130	7.14633621765187e-06\\
131	6.89964623558094e-06\\
132	6.64150874185355e-06\\
133	6.28968287569217e-06\\
134	6.12965068301528e-06\\
135	5.96887100119845e-06\\
136	5.60720954444479e-06\\
137	5.4029447150641e-06\\
138	5.26612693971471e-06\\
139	5.14107365370069e-06\\
140	4.96479143681058e-06\\
141	4.80648247535e-06\\
142	4.54997130250425e-06\\
143	4.31485011775576e-06\\
144	4.25525664509935e-06\\
145	4.05983632165914e-06\\
146	3.86038639069763e-06\\
147	3.76999760648614e-06\\
148	3.68584352982437e-06\\
149	3.57966673476537e-06\\
150	3.40413335856154e-06\\
151	3.25207296662623e-06\\
152	3.06208800509487e-06\\
153	3.00832824750076e-06\\
154	2.89737766783253e-06\\
155	2.82838426574094e-06\\
156	2.71874324204139e-06\\
157	2.61723955634396e-06\\
158	2.3993224962735e-06\\
159	2.35523933116239e-06\\
160	2.31633151351664e-06\\
161	2.22009279633956e-06\\
162	2.14864976152266e-06\\
163	2.03259251486615e-06\\
164	1.9446011690243e-06\\
165	1.87844625201059e-06\\
166	1.83198277795747e-06\\
167	1.80028865615432e-06\\
168	1.71658553275752e-06\\
169	1.66238852875046e-06\\
170	1.60798314465366e-06\\
171	1.56109012010876e-06\\
172	1.49356535297194e-06\\
173	1.45812348506236e-06\\
174	1.40925257792599e-06\\
175	1.32262914365137e-06\\
176	1.31310306796351e-06\\
177	1.28522623478686e-06\\
178	1.20450214737034e-06\\
179	1.16108616205869e-06\\
180	1.12517557240663e-06\\
181	1.06287014164787e-06\\
182	1.0163791937677e-06\\
183	9.90573094536691e-07\\
184	9.8027334693373e-07\\
185	9.17060327779605e-07\\
186	8.90435253238347e-07\\
187	8.53057386949294e-07\\
188	8.33702444585276e-07\\
189	8.14826610670657e-07\\
190	7.69992760126944e-07\\
191	7.27784044175201e-07\\
192	7.05562268356147e-07\\
193	6.85963854652734e-07\\
194	6.43314942611522e-07\\
195	6.41804313356991e-07\\
196	6.04439051104562e-07\\
197	5.9931697868151e-07\\
198	5.87918408384794e-07\\
199	5.43073213711911e-07\\
200	5.11990171963662e-07\\
201	4.99903776963653e-07\\
202	4.86928759870241e-07\\
203	4.82913032432246e-07\\
204	4.62893122402631e-07\\
205	4.37662435342423e-07\\
206	4.02303129796808e-07\\
207	3.97516105367399e-07\\
208	3.87419310710304e-07\\
209	3.8468011404169e-07\\
210	3.55195057706422e-07\\
211	3.40577896742576e-07\\
212	3.30864959387372e-07\\
213	3.2693878872503e-07\\
214	3.10115154030935e-07\\
215	2.95511296774291e-07\\
216	2.87897407352654e-07\\
217	2.72171174330147e-07\\
218	2.69239201321003e-07\\
219	2.60264878621382e-07\\
220	2.42832792324577e-07\\
221	2.37333277126402e-07\\
222	2.22017104803233e-07\\
223	2.18243844936952e-07\\
224	2.04796747793735e-07\\
225	1.99428713857822e-07\\
226	1.93120597498487e-07\\
227	1.85013266623081e-07\\
228	1.80738060630428e-07\\
229	1.75727810002104e-07\\
230	1.6795288766253e-07\\
231	1.58686467316607e-07\\
232	1.55483337616192e-07\\
233	1.48297158874414e-07\\
234	1.44111667515087e-07\\
235	1.36902166379153e-07\\
236	1.32224872801407e-07\\
237	1.27028133886924e-07\\
238	1.21774459851733e-07\\
239	1.1921144336344e-07\\
240	1.09838079059692e-07\\
241	1.07096400581624e-07\\
242	1.02038033183505e-07\\
243	9.87687839080493e-08\\
244	9.78800565730585e-08\\
245	9.05804238809261e-08\\
246	8.97419453690661e-08\\
247	8.50149414234051e-08\\
248	8.26887591019273e-08\\
249	7.92836308013051e-08\\
250	7.82032542724592e-08\\
251	7.41441472622351e-08\\
252	6.92182194072483e-08\\
253	6.66834122576936e-08\\
254	6.50179480868551e-08\\
255	6.25913574606846e-08\\
256	6.18932890120222e-08\\
257	5.7276378983586e-08\\
258	5.61343803946475e-08\\
259	5.38700514046381e-08\\
260	5.13051549047378e-08\\
261	4.78044097851069e-08\\
262	4.71940098996617e-08\\
263	4.61305962321376e-08\\
264	4.47129791743456e-08\\
265	4.12139727544857e-08\\
266	3.93170698140513e-08\\
267	3.858209664464e-08\\
268	3.67053493361452e-08\\
269	3.63415235877885e-08\\
270	3.35782181249288e-08\\
271	3.14581688275855e-08\\
272	3.12275391249867e-08\\
273	2.98803379380755e-08\\
274	2.91511078874648e-08\\
275	2.74084801409681e-08\\
276	2.71479579418641e-08\\
277	2.56133759377574e-08\\
278	2.45350847763197e-08\\
279	2.44932819590896e-08\\
280	2.38676263625577e-08\\
281	2.1910239112238e-08\\
282	2.08883780792739e-08\\
283	2.04727677266281e-08\\
284	1.97195252541789e-08\\
285	1.89517779220232e-08\\
286	1.75484139046369e-08\\
287	1.74797741009403e-08\\
288	1.65049318207373e-08\\
289	1.61080433817119e-08\\
290	1.49955967800231e-08\\
291	1.44717503638531e-08\\
292	1.39232692858287e-08\\
293	1.33140391190463e-08\\
294	1.31132774502853e-08\\
295	1.19994155121737e-08\\
296	1.17080422428049e-08\\
297	1.1109094607268e-08\\
298	1.09640860272682e-08\\
299	1.05482151467425e-08\\
300	9.466266590155e-09\\
};
\end{axis}
\end{tikzpicture}%

%% file: fig/Pex1/rel_err_training.tikz
\definecolor{mycolor1}{rgb}{0.00000,0.44700,0.74100}%
\definecolor{mycolor2}{rgb}{0.85000,0.32500,0.09800}%
\begin{tikzpicture}

\begin{axis}[%
width=\fwidth,
height=0.999\fheight,
at={(0\fwidth,0\fheight)},
scale only axis,
xmin=25,
xmax=75,
xlabel style={font=\color{white!15!black}},
xlabel={Reduced dimension},
ymode=log,
ymin=0.0001,
ymax=0.01,
yminorticks=true,
ylabel style={font=\color{white!15!black}},
ylabel={Relative error},
axis background/.style={fill=white},
xmajorgrids,
xminorgrids,
ymajorgrids,
yminorgrids,
legend style={legend cell align=left, align=left, draw=white!15!black}
]
\addplot [color=mycolor1, line width=2.0pt,dashed, mark size=2.5pt, mark=triangle, mark repeat=10, mark options={solid, mycolor1}]
  table[row sep=crcr]{%
25	0.00278639814217809\\
26	0.0026045666509328\\
27	0.00249013405232446\\
28	0.0023097037274231\\
29	0.00216774248199421\\
30	0.00202349173157815\\
31	0.00190889315909088\\
32	0.00261437639893268\\
33	0.00169496186442185\\
34	0.00160449965239008\\
35	0.00152507147048109\\
36	0.00144554685634393\\
37	0.0013671285960693\\
38	0.00130032382812657\\
39	0.0012305297035218\\
40	0.00117033830076947\\
41	0.00110969865405596\\
42	0.00104820934951787\\
43	0.000992901008586423\\
44	0.000945662011769559\\
45	0.000900849053978366\\
46	0.000855507768591343\\
47	0.000814957978430928\\
48	0.000774248004877222\\
49	0.000735539413052992\\
50	0.000699294715102331\\
51	0.000666686235929648\\
52	0.000634810571143711\\
53	0.000599695154614958\\
54	0.000568998775627424\\
55	0.000542901581543289\\
56	0.000515987333229457\\
57	0.000490627906702614\\
58	0.000467244419214519\\
59	0.000444676436787724\\
60	0.000429935965555538\\
61	0.000417652981780974\\
62	0.000404998842222395\\
63	0.000393069537938404\\
64	0.000375996206882611\\
65	0.000343365518529157\\
66	0.000344133125473528\\
67	0.00032446491911601\\
68	0.000300862138162209\\
69	0.000314340235798676\\
70	0.000276861097514574\\
71	0.000281677874848137\\
72	0.000267504310948461\\
73	0.000263333551597236\\
74	0.000251747296281611\\
75	0.000332159776528106\\
};
\addlegendentry{Inference tQR}

\addplot [color=mycolor2, line width=2.0pt, dotted, mark size=2.5pt, mark=o, mark repeat=10, mark phase = 5, mark options={solid, mycolor2}]
  table[row sep=crcr]{%
25	0.00278625606891295\\
26	0.00260449793455654\\
27	0.00263395711100096\\
28	0.00230945373247215\\
29	0.00216754767987239\\
30	0.00202338800579541\\
31	0.00190884941769992\\
32	0.00676861163580303\\
33	0.00169462267087534\\
34	0.00160404998804985\\
35	0.00152459548696303\\
36	0.00144475013914247\\
37	0.00136668985001392\\
38	0.00129964949069938\\
39	0.00122987185927761\\
40	0.00116963616780193\\
41	0.00110891478475615\\
42	0.00104724059571846\\
43	0.000992206774940555\\
44	0.000945054150115406\\
45	0.000898945590086643\\
46	0.000853908807150122\\
47	0.000813163656412565\\
48	0.000772718206850338\\
49	0.000733592078693282\\
50	0.000697554343465087\\
51	0.000663740750214752\\
52	0.000630299246904725\\
53	0.000597103174422143\\
54	0.000567184564966592\\
55	0.000539395704143923\\
56	0.000510983194928409\\
57	0.000485559022777611\\
58	0.000462209332420429\\
59	0.00043873057235489\\
60	0.000418482237575633\\
61	0.000398585518783334\\
62	0.00037967320106869\\
63	0.000361891485334449\\
64	0.000345373209744564\\
65	0.000330006566104594\\
66	0.000315968153896597\\
67	0.000301857642662125\\
68	0.000288268256268311\\
69	0.000275497884045517\\
70	0.000263475114951635\\
71	0.000251752651672039\\
72	0.000240947278497052\\
73	0.000230269330347861\\
74	0.000219714324779092\\
75	0.000209799497222883\\
};
\addlegendentry{Inference Tikhonov}

\end{axis}
\end{tikzpicture}%

%% file: fig/Pex1/param_err.tikz
\begin{tikzpicture}

\begin{axis}[%
width=\fwidth,
height=0.999\fheight,
at={(0\fwidth,0\fheight)},
scale only axis,
xmin=0.5,
xmax=1.1,
xlabel style={font=\color{white!15!black}},
xlabel={Parameters set},
ymode=log,
ymin=1e-05,
ymax=0.1,
yminorticks=true,
ylabel style={font=\color{white!15!black}},
ylabel={Relative error},
axis background/.style={fill=white},
xmajorgrids,
ymajorgrids,
yminorgrids,
legend style={legend cell align=left, align=left, draw=white!15!black}
]
\addplot [color=black, line width=2.0pt, draw=none, only marks, mark size=2.5pt, mark=o, mark options={solid, black}]
  table[row sep=crcr]{%
0.589048622548086	0.00264982447136121\\
0.719948316447661	0.00237363305829157\\
0.850848010347236	0.00223115959065009\\
0.98174770424681	0.0025350293601419\\
};
\addlegendentry{Test Parameters}

\addplot [color=black, line width=2.0pt, draw=none, only marks, mark size=1.8pt, mark=square, mark options={solid, black}]
  table[row sep=crcr]{%
0.523598775598299	0.00366313368881339\\
0.654498469497873	0.00245709873803542\\
0.785398163397448	0.00227634944827331\\
0.916297857297023	0.00229876368509079\\
1.0471975511966	0.00297906772278355\\
};
\addlegendentry{Training Parameters}

\addplot [color=red, line width=2.0pt, draw=none, only marks, mark size=2.5pt, mark=o, mark options={solid, red}, forget plot]
  table[row sep=crcr]{%
0.589048622548086	0.000483957378892707\\
0.719948316447661	0.000272754561191688\\
0.850848010347236	0.000212674405441802\\
0.98174770424681	0.000192732136760122\\
};
\addplot [color=red, line width=2.0pt, draw=none, only marks, mark size=1.8pt, mark=square, mark options={solid, red}, forget plot]
  table[row sep=crcr]{%
0.523598775598299	0.000518863667702133\\
0.654498469497873	0.000365691518580049\\
0.785398163397448	0.000224654193839226\\
0.916297857297023	0.000212628109703753\\
1.0471975511966	0.000228075450384183\\
};
\end{axis}
\end{tikzpicture}%

%% file: fig/Nex2/rel_err.tikz
\definecolor{mycolor1}{rgb}{0.00000,0.44700,0.74100}%
\definecolor{mycolor2}{rgb}{0.85000,0.32500,0.09800}%
\definecolor{mycolor3}{rgb}{0.92900,0.69400,0.12500}%
\begin{tikzpicture}

\begin{axis}[%
width=0.995\fwidth,
height=\fheight,
at={(0\fwidth,0\fheight)},
scale only axis,
xmin=15,
xmax=50,
xlabel style={font=\color{white!15!black}},
xlabel={Reduced dimension},
ymode=log,
ymin=0.00001,
ymax=0.01,
yminorticks=true,
ylabel style={font=\color{white!15!black}},
ylabel={Relative error},
axis background/.style={fill=white},
xmajorgrids,
ymajorgrids,
yminorgrids,
legend style={legend cell align=left, align=left, draw=white!15!black}
]
\addplot [color=mycolor1, line width=2.0pt]
  table[row sep=crcr]{%
15	0.00164146288165366\\
16	0.00231454608797397\\
17	0.00454923186017515\\
18	0.00203824182606238\\
19	0.00254942935329309\\
20	0.00133556630810592\\
21	0.00125002097689712\\
22	0.000979842546741085\\
23	0.00106941317942511\\
24	0.00102052637791345\\
25	0.000577042243392006\\
26	0.000914019120898893\\
27	0.000472594311812238\\
28	0.000468750939761047\\
29	0.000463203654165794\\
30	0.000921264648192169\\
31	0.000910307787005673\\
32	0.000953875660810869\\
33	0.00060883891666413\\
34	0.00118207488282618\\
35	0.000608552239063065\\
36	0.0010532439522782\\
37	0.000611991663230374\\
38	0.000607133270342172\\
39	0.000444901945631638\\
40	0.000199664229942492\\
41	0.000120378586249188\\
42	0.000122618718502065\\
43	9.97792076224462e-05\\
44	9.58416013022434e-05\\
45	8.88361658625882e-05\\
46	8.58772833921628e-05\\
47	6.72201776251672e-05\\
48	7.56141162433707e-05\\
49	6.15186969842296e-05\\
50	5.54200719848054e-05\\
};
\addlegendentry{POD}

\addplot [color=mycolor2, line width=2.0pt,dashed, mark size=2.0pt, mark=triangle, mark repeat=10, mark options={solid, mycolor2}]
  table[row sep=crcr]{%
15	0.00152094140838023\\
16	0.00139411999414811\\
17	0.00126727411679246\\
18	0.00113247345672082\\
19	0.000989970315519913\\
20	0.000848778614437691\\
21	0.000765266020189398\\
22	0.000675144985838956\\
23	0.000581984069930702\\
24	0.000501211293878488\\
25	0.000437016292080964\\
26	0.000401094180336563\\
27	0.000369369246202197\\
28	0.000347773479888117\\
29	0.000325164076667048\\
30	0.000302152570793941\\
31	0.000281073008916211\\
32	0.000259832713968677\\
33	0.000237973089093083\\
34	0.0002150168134215\\
35	0.000192068599058465\\
36	0.000167802240759559\\
37	0.000144917819954246\\
38	0.000131392244567537\\
39	0.00011928506761866\\
40	0.000107805816898617\\
41	9.77685898032673e-05\\
42	8.90621814926593e-05\\
43	7.98113117797363e-05\\
44	7.13456941753774e-05\\
45	6.30744509415308e-05\\
46	5.71113786105031e-05\\
47	5.19577660081306e-05\\
48	4.66274811690371e-05\\
49	4.10458146842377e-05\\
50	3.5792764888218e-05\\
};
\addlegendentry{Inference tQR}

\addplot [color=mycolor3, line width=2.0pt, dotted, mark size=2.0pt, mark=o, mark repeat=10, mark phase = 5, mark options={solid, mycolor3}]
  table[row sep=crcr]{%
15	0.00152093274971325\\
16	0.00139372763396822\\
17	0.00126725381569234\\
18	0.00113246849602611\\
19	0.000989965428402032\\
20	0.000848774013591304\\
21	0.000765261494464139\\
22	0.000675141437559292\\
23	0.000581977668128481\\
24	0.000501205093613601\\
25	0.000437010228642799\\
26	0.000401085516820253\\
27	0.000369357694453227\\
28	0.000347733907054426\\
29	0.000325147988864735\\
30	0.00030213778028046\\
31	0.000280998157322982\\
32	0.000259784575466767\\
33	0.000237966410172106\\
34	0.000214994892620334\\
35	0.000192059614390317\\
36	0.000167788075577908\\
37	0.000144898087447383\\
38	0.000131376513974997\\
39	0.000119245650520284\\
40	0.00010777659505875\\
41	9.77096848318744e-05\\
42	8.90133312312077e-05\\
43	7.97554687224026e-05\\
44	7.12438665964135e-05\\
45	6.29291998156636e-05\\
46	5.69582421202213e-05\\
47	5.16963623853582e-05\\
48	4.63558227893282e-05\\
49	4.06405000162602e-05\\
50	3.51694395616465e-05\\
};
\addlegendentry{Inference Tikhonov}

\end{axis}
\end{tikzpicture}%

%% file: fig/Nex2/sing.tikz
\definecolor{mycolor1}{rgb}{0.00000,0.44700,0.74100}%
\begin{tikzpicture}

\begin{axis}[%
width=0.961\fwidth,
height=\fheight,
at={(0\fwidth,0\fheight)},
scale only axis,
xmin=1,
xmax=300,
xlabel style={font=\color{white!15!black}},
xlabel={$k$},
ymode=log,
ymin=2.69054960474281e-11,
ymax=1,
yminorticks=true,
ylabel style={font=\color{white!15!black}},
ylabel={$\sigma_k / \sigma_1$},
axis background/.style={fill=white},
xmajorgrids,
xminorgrids,
ymajorgrids,
yminorgrids
]
\addplot [color=mycolor1, line width=2.0pt, forget plot]
  table[row sep=crcr]{%
1	1\\
2	0.0370303867331579\\
3	0.0359732140792997\\
4	0.021569843459679\\
5	0.00657820672938038\\
6	0.00650937512526578\\
7	0.0029135511250204\\
8	0.00173697739843632\\
9	0.00168265412826634\\
10	0.00108483331859857\\
11	0.00107148537319668\\
12	0.000998014160436999\\
13	0.000974561669655522\\
14	0.000899873802692618\\
15	0.000707218985200544\\
16	0.000609953177341252\\
17	0.000580993600921846\\
18	0.000569642667719892\\
19	0.000550847231589929\\
20	0.000510347851458616\\
21	0.000367734579878256\\
22	0.000360873386050639\\
23	0.000342778711002485\\
24	0.000296268046855744\\
25	0.000245812141470012\\
26	0.000173798397776987\\
27	0.000156600331909381\\
28	0.000124726313906366\\
29	0.000123481683575843\\
30	0.000120332068265394\\
31	0.000111208246979067\\
32	0.000107291924341445\\
33	0.000104421245904801\\
34	0.000102243273794383\\
35	9.66439895157714e-05\\
36	9.3607440579089e-05\\
37	8.47382418939906e-05\\
38	6.12191242697995e-05\\
39	5.52281677280189e-05\\
40	5.11094726080903e-05\\
41	4.55550004097266e-05\\
42	4.03615401423627e-05\\
43	3.95924205032058e-05\\
44	3.59079843891061e-05\\
45	3.34554787780005e-05\\
46	2.69288725413861e-05\\
47	2.38017538553481e-05\\
48	2.29214806039811e-05\\
49	2.23345282716565e-05\\
50	2.03995052216721e-05\\
51	1.62693558475444e-05\\
52	1.38736855079689e-05\\
53	1.35379962646603e-05\\
54	1.16893712083828e-05\\
55	1.01943871778077e-05\\
56	9.46879100151887e-06\\
57	7.5055069249839e-06\\
58	6.13600949427472e-06\\
59	5.87288560728397e-06\\
60	5.80008653282717e-06\\
61	5.3162280385361e-06\\
62	4.16421939014786e-06\\
63	3.57773660849991e-06\\
64	3.46106927723182e-06\\
65	2.95158318986842e-06\\
66	2.38069568490034e-06\\
67	2.22119464358453e-06\\
68	2.00915879534269e-06\\
69	1.60665272471225e-06\\
70	1.59139302885079e-06\\
71	1.5156134861689e-06\\
72	1.14924536346685e-06\\
73	9.65874736964018e-07\\
74	9.06727931527111e-07\\
75	9.02912305727913e-07\\
76	7.83529476909614e-07\\
77	5.87621995194117e-07\\
78	5.02284124409801e-07\\
79	4.7715985886453e-07\\
80	4.16879422629449e-07\\
81	4.05440699589974e-07\\
82	3.54344243652136e-07\\
83	3.16235541912133e-07\\
84	3.15403090257307e-07\\
85	2.8527513319522e-07\\
86	2.75319177537701e-07\\
87	2.46636425235773e-07\\
88	2.13393532644331e-07\\
89	1.96935712152697e-07\\
90	1.96391483314322e-07\\
91	1.73433691673239e-07\\
92	1.56646439433639e-07\\
93	1.56251839014807e-07\\
94	1.41166984421925e-07\\
95	1.39715602412269e-07\\
96	1.30122298838231e-07\\
97	1.25014691223457e-07\\
98	1.24720620250943e-07\\
99	1.13886856417223e-07\\
100	1.13797298448843e-07\\
101	1.0207310153321e-07\\
102	1.01973053900832e-07\\
103	1.00479214338607e-07\\
104	9.1671075353175e-08\\
105	8.80146727051867e-08\\
106	8.79466111950834e-08\\
107	8.0846466354788e-08\\
108	7.86257924915363e-08\\
109	7.83624302602116e-08\\
110	6.95291176642989e-08\\
111	6.84951736179343e-08\\
112	6.69670223373995e-08\\
113	6.6654583421706e-08\\
114	6.35840733583256e-08\\
115	6.19243745378974e-08\\
116	6.11769619050812e-08\\
117	5.07918928812249e-08\\
118	4.75727183376314e-08\\
119	4.74722975468713e-08\\
120	4.30484183191703e-08\\
121	3.97464485849423e-08\\
122	3.8860270863609e-08\\
123	3.76280880604983e-08\\
124	3.69973660873594e-08\\
125	3.34570938186296e-08\\
126	3.34546952087972e-08\\
127	3.17292523965592e-08\\
128	3.09590701604692e-08\\
129	3.09196033321571e-08\\
130	2.57647023453973e-08\\
131	2.49774249505096e-08\\
132	2.4944554933938e-08\\
133	2.40717958489636e-08\\
134	2.39896939045146e-08\\
135	2.38325582849768e-08\\
136	2.23346330398406e-08\\
137	2.22265846907525e-08\\
138	2.14468340303142e-08\\
139	1.99600141935665e-08\\
140	1.83913044450213e-08\\
141	1.79779426237586e-08\\
142	1.78692171554675e-08\\
143	1.63980110367724e-08\\
144	1.63218459262733e-08\\
145	1.56904363588311e-08\\
146	1.42731688741248e-08\\
147	1.40839252110878e-08\\
148	1.26109711892324e-08\\
149	1.23927336801113e-08\\
150	1.20412924382113e-08\\
151	1.10492888991902e-08\\
152	1.08339512540161e-08\\
153	1.081287640772e-08\\
154	1.05779785030445e-08\\
155	9.1070680757609e-09\\
156	8.11252070969119e-09\\
157	7.98508126931516e-09\\
158	7.90143072510081e-09\\
159	7.37431753871366e-09\\
160	7.24336131729933e-09\\
161	7.12067122087502e-09\\
162	6.5916890936648e-09\\
163	6.08329759018448e-09\\
164	6.04089714384461e-09\\
165	5.76316813342416e-09\\
166	5.22619234597003e-09\\
167	5.16195756954489e-09\\
168	5.02791145073873e-09\\
169	4.49729012906609e-09\\
170	4.32078350974581e-09\\
171	4.28423373543995e-09\\
172	3.84138241984077e-09\\
173	3.54148510992301e-09\\
174	3.52395153122717e-09\\
175	3.48659949805131e-09\\
176	3.3517825392196e-09\\
177	3.11408462469773e-09\\
178	3.02201379252311e-09\\
179	2.69572257150955e-09\\
180	2.46426015811701e-09\\
181	2.43187379797602e-09\\
182	2.26952627126767e-09\\
183	2.06614788612782e-09\\
184	2.05501949403136e-09\\
185	2.03315822896448e-09\\
186	1.8482928918696e-09\\
187	1.78596432697907e-09\\
188	1.76234756979263e-09\\
189	1.66375142222005e-09\\
190	1.64972643389505e-09\\
191	1.63596664310606e-09\\
192	1.46084562173314e-09\\
193	1.4353126906316e-09\\
194	1.42091978061925e-09\\
195	1.29016917128152e-09\\
196	1.26879970911606e-09\\
197	1.2275280329646e-09\\
198	1.18112703755734e-09\\
199	1.04506989922699e-09\\
200	1.03746205750178e-09\\
201	9.82887894436909e-10\\
202	9.10679262247483e-10\\
203	9.01574786765771e-10\\
204	8.54653242460492e-10\\
205	7.73359678233584e-10\\
206	7.29167665381007e-10\\
207	7.22455237392301e-10\\
208	6.88637000180485e-10\\
209	6.64766789153035e-10\\
210	6.5955193045572e-10\\
211	6.11973541204316e-10\\
212	5.7625415102238e-10\\
213	5.57866234445709e-10\\
214	5.44141340200376e-10\\
215	5.04595705433039e-10\\
216	5.03825972423778e-10\\
217	4.88715304930799e-10\\
218	4.55269893871353e-10\\
219	4.48111680102738e-10\\
220	4.36430109743898e-10\\
221	3.80462991563038e-10\\
222	3.57798601610826e-10\\
223	3.47844460916714e-10\\
224	3.34870508317563e-10\\
225	3.16573302958634e-10\\
226	3.09154759184788e-10\\
227	3.06105036684999e-10\\
228	2.75520783249573e-10\\
229	2.73612366168897e-10\\
230	2.59949399020542e-10\\
231	2.36892579081201e-10\\
232	2.34233555786213e-10\\
233	2.19959505247945e-10\\
234	2.03495439240697e-10\\
235	1.95750935225638e-10\\
236	1.90790241169057e-10\\
237	1.76073148242707e-10\\
238	1.69544939869957e-10\\
239	1.65267508442259e-10\\
240	1.5589989115599e-10\\
241	1.53506674656945e-10\\
242	1.5169598135179e-10\\
243	1.38485500594609e-10\\
244	1.31773344170711e-10\\
245	1.30284188423762e-10\\
246	1.25755250784344e-10\\
247	1.15003066399348e-10\\
248	1.075327734394e-10\\
249	1.05844061954463e-10\\
250	9.75105531753512e-11\\
251	9.43744469725142e-11\\
252	9.28424965829801e-11\\
253	8.65117386052806e-11\\
254	8.48064033511329e-11\\
255	8.4291727640273e-11\\
256	8.30277766790626e-11\\
257	8.1094550598871e-11\\
258	7.7703468247639e-11\\
259	7.60936807703556e-11\\
260	7.3943575147735e-11\\
261	7.07101574939772e-11\\
262	6.98336657415926e-11\\
263	6.81808700768911e-11\\
264	6.77833321408853e-11\\
265	6.63401976243083e-11\\
266	6.29245729444987e-11\\
267	6.15773515277849e-11\\
268	6.04191060119078e-11\\
269	5.75798052398248e-11\\
270	5.5133203636508e-11\\
271	5.43052432095484e-11\\
272	5.30684902952731e-11\\
273	5.25910232791044e-11\\
274	5.21711092210013e-11\\
275	5.09040070160809e-11\\
276	4.96496879615423e-11\\
277	4.89478369877162e-11\\
278	4.72085078259792e-11\\
279	4.52585879928883e-11\\
280	4.45352991381056e-11\\
281	4.25595941593727e-11\\
282	4.2446287857869e-11\\
283	4.10961846834174e-11\\
284	4.08910033807336e-11\\
285	3.92988038201601e-11\\
286	3.91809434020907e-11\\
287	3.78797997913962e-11\\
288	3.68755956204489e-11\\
289	3.6553333047237e-11\\
290	3.63038204828224e-11\\
291	3.53376248549064e-11\\
292	3.45913230277203e-11\\
293	3.31770320166276e-11\\
294	3.27611106552332e-11\\
295	3.17647939387073e-11\\
296	2.97344327133876e-11\\
297	2.89524298729931e-11\\
298	2.77451814769891e-11\\
299	2.74913133386712e-11\\
300	2.69054960474281e-11\\
};
\end{axis}
\end{tikzpicture}%

%% file: fig/Pex2/rel_err_training.tikz
\definecolor{mycolor1}{rgb}{0.00000,0.44700,0.74100}%
\definecolor{mycolor2}{rgb}{0.85000,0.32500,0.09800}%
\begin{tikzpicture}

\begin{axis}[%
width=\fwidth,
height=0.999\fheight,
at={(0\fwidth,0\fheight)},
scale only axis,
xmin=20,
xmax=65,
xlabel style={font=\color{white!15!black}},
xlabel={Reduced dimension},
ymode=log,
ymin=5e-05,
ymax=0.01,
yminorticks=true,
ylabel style={font=\color{white!15!black}},
ylabel={Relative error},
axis background/.style={fill=white},
xmajorgrids,
xminorgrids,
ymajorgrids,
yminorgrids,
legend style={legend cell align=left, align=left, draw=white!15!black}
]
\addplot [color=mycolor1, line width=2.0pt,dashed, mark size=2.5pt, mark=triangle, mark repeat=10, mark options={solid, mycolor1}]
  table[row sep=crcr]{%
20	0.000879585678974916\\
21	0.000827256466673637\\
22	0.000747674563707237\\
23	0.000686832462958634\\
24	0.00063314192604609\\
25	0.000590565483386273\\
26	0.00054790065901951\\
27	0.000514934907039756\\
28	0.000482986388251947\\
29	0.000450533987549361\\
30	0.000419790002915325\\
31	0.000392699322837922\\
32	0.000365488339363763\\
33	0.000342365691521137\\
34	0.000320141924580333\\
35	0.000301448325920208\\
36	0.00028369219730475\\
37	0.000269142138356234\\
38	0.000253736403755828\\
39	0.000239781076088906\\
40	0.000227263862901586\\
41	0.000216562311025758\\
42	0.000206233482738893\\
43	0.000195670721955695\\
44	0.000185878677542517\\
45	0.000178732530228323\\
46	0.000214555178938913\\
47	0.000163187622894893\\
48	0.000161418171466337\\
49	0.000192678127021278\\
50	0.000144683034957852\\
51	0.000134136563975173\\
52	0.000129545599418903\\
53	0.000126175077709557\\
54	0.000121842178034369\\
55	0.000116502348359421\\
56	0.000112374440899026\\
57	0.000105498738871593\\
58	9.88783186062935e-05\\
59	0.000100762539142464\\
60	9.40076487670722e-05\\
61	8.66356149662286e-05\\
62	8.09496966339987e-05\\
63	8.31203997927011e-05\\
64	7.42466011669585e-05\\
65	7.41849999900745e-05\\
};
\addlegendentry{Inference tQR}

\addplot [color=mycolor2, line width=2.0pt, dotted, mark size=2.5pt, mark=o, mark repeat=10, mark phase = 5, mark options={solid, mycolor2}]
  table[row sep=crcr]{%
20	0.000877241169688815\\
21	0.000806532701744536\\
22	0.000745297218982737\\
23	0.000685451076473111\\
24	0.000630519297408053\\
25	0.000585635257443501\\
26	0.000545970317453055\\
27	0.000511971213923486\\
28	0.000480023167654034\\
29	0.000447944938179314\\
30	0.000419229630419145\\
31	0.000392254952557699\\
32	0.000365081686828177\\
33	0.000341632998978316\\
34	0.000318633910942116\\
35	0.000300968893602373\\
36	0.000282798116258903\\
37	0.000268007718418338\\
38	0.000252615199814123\\
39	0.000238865145782727\\
40	0.000225710853653885\\
41	0.000214032238430591\\
42	0.000202703740100432\\
43	0.000192814735688813\\
44	0.000183359264882344\\
45	0.000173760041022705\\
46	0.000164629563089928\\
47	0.00015699393724286\\
48	0.000149252209969644\\
49	0.000142001327707275\\
50	0.000134785829402445\\
51	0.000127828612619754\\
52	0.000121743256334423\\
53	0.000115271707114557\\
54	0.000110719414305181\\
55	0.000103829825186751\\
56	9.82557173926434e-05\\
57	9.34360682957958e-05\\
58	8.84693011696665e-05\\
59	8.38927299729405e-05\\
60	7.92841939418971e-05\\
61	7.50369092615576e-05\\
62	7.08286046428677e-05\\
63	6.73374172390842e-05\\
64	6.42642064026673e-05\\
65	6.12092386611556e-05\\
};
\addlegendentry{Inference Tikhonov}

\end{axis}
\end{tikzpicture}%

%% file: fig/Pex2/sing.tikz
\definecolor{mycolor1}{rgb}{0.00000,0.44700,0.74100}%
\begin{tikzpicture}

\begin{axis}[%
width=0.961\fwidth,
height=\fheight,
at={(0\fwidth,0\fheight)},
scale only axis,
xmin=0,
xmax=300,
xlabel style={font=\color{white!15!black}},
xlabel={$k$},
ymode=log,
ymin=1e-10,
ymax=1,
yminorticks=true,
ylabel style={font=\color{white!15!black}},
ylabel={$\sigma_k / \sigma_1$},
axis background/.style={fill=white},
xmajorgrids,
xminorgrids,
ymajorgrids,
yminorgrids
]
\addplot [color=mycolor1, line width=2.0pt, forget plot]
  table[row sep=crcr]{%
1	1\\
2	0.038555630708281\\
3	0.0365977453055534\\
4	0.0275891956851907\\
5	0.0149255921213209\\
6	0.00711354192471883\\
7	0.00701579720193136\\
8	0.00256508518961198\\
9	0.00217834561165005\\
10	0.0014240568790174\\
11	0.00131066209663764\\
12	0.00130220994131421\\
13	0.0011726071357654\\
14	0.00111568032945697\\
15	0.00101105828660569\\
16	0.000914043286727094\\
17	0.000702288236314229\\
18	0.00068089392895881\\
19	0.000646044683868272\\
20	0.000546155640879275\\
21	0.000345556399773151\\
22	0.000308726335169667\\
23	0.000293167547822867\\
24	0.00026938819708272\\
25	0.000234099333065826\\
26	0.000212269096703009\\
27	0.000190099468216272\\
28	0.000178298758982017\\
29	0.000172867394963603\\
30	0.000158092069352924\\
31	0.000148241972623324\\
32	0.000143737783088101\\
33	0.000128982795531012\\
34	0.00012347096032171\\
35	0.000104825686451125\\
36	0.000103193979839104\\
37	9.04372338504951e-05\\
38	8.96969234456831e-05\\
39	8.23665631127351e-05\\
40	7.83308026502844e-05\\
41	7.18024011693563e-05\\
42	6.88480094771484e-05\\
43	6.26657808358985e-05\\
44	5.97593631782577e-05\\
45	5.86501957872812e-05\\
46	5.58319961457307e-05\\
47	4.98553307475684e-05\\
48	4.87917419524267e-05\\
49	4.63698295024207e-05\\
50	4.45678209467532e-05\\
51	4.22952213873069e-05\\
52	3.94245975427323e-05\\
53	3.92082619304027e-05\\
54	3.57447982928342e-05\\
55	3.48854686004751e-05\\
56	3.36208558916603e-05\\
57	3.04561405282811e-05\\
58	3.01118257052769e-05\\
59	2.80552742585793e-05\\
60	2.74882916684889e-05\\
61	2.56840494172638e-05\\
62	2.48298685257502e-05\\
63	2.2008356411772e-05\\
64	2.01326387597518e-05\\
65	1.96154059295214e-05\\
66	1.87318882680707e-05\\
67	1.81326397535955e-05\\
68	1.66062612030341e-05\\
69	1.65373309873988e-05\\
70	1.44919221972594e-05\\
71	1.43845016205007e-05\\
72	1.35972049735955e-05\\
73	1.1981467086461e-05\\
74	1.1607427576453e-05\\
75	1.13898234492386e-05\\
76	1.0760241979237e-05\\
77	1.06144795143982e-05\\
78	1.03005762837517e-05\\
79	9.69909981020851e-06\\
80	9.4686319458391e-06\\
81	9.36381791772237e-06\\
82	8.9190191627029e-06\\
83	8.7545255629109e-06\\
84	7.78973433019442e-06\\
85	7.3471476040109e-06\\
86	7.24883920276568e-06\\
87	7.0201689824183e-06\\
88	6.49709138445209e-06\\
89	6.26409227991121e-06\\
90	6.09874525937215e-06\\
91	5.9632618323966e-06\\
92	5.79276260783679e-06\\
93	5.37829419664067e-06\\
94	5.04091835934158e-06\\
95	4.84926753086342e-06\\
96	4.30504200809581e-06\\
97	4.24393031464631e-06\\
98	4.07596933316825e-06\\
99	3.99147378371979e-06\\
100	3.90005820178191e-06\\
101	3.58487713839644e-06\\
102	3.52668591007057e-06\\
103	3.47559361003373e-06\\
104	3.21792600179427e-06\\
105	3.18153073239493e-06\\
106	3.02895567137365e-06\\
107	3.00016991218393e-06\\
108	2.88208083060025e-06\\
109	2.76061808127458e-06\\
110	2.66363705630276e-06\\
111	2.58798079952421e-06\\
112	2.52060940100751e-06\\
113	2.36107086618748e-06\\
114	2.30990095143132e-06\\
115	2.29037424073586e-06\\
116	2.19723454955543e-06\\
117	2.09158757907936e-06\\
118	1.94602953992209e-06\\
119	1.92484374351668e-06\\
120	1.87677206957707e-06\\
121	1.76223288833381e-06\\
122	1.68401352496979e-06\\
123	1.63218564651263e-06\\
124	1.58155424324628e-06\\
125	1.46970508427959e-06\\
126	1.41975257775743e-06\\
127	1.36895870562667e-06\\
128	1.33375685084003e-06\\
129	1.30904767095617e-06\\
130	1.26751822458808e-06\\
131	1.17558759468142e-06\\
132	1.09823574031619e-06\\
133	1.08436237106971e-06\\
134	1.04962571847747e-06\\
135	1.0249233917465e-06\\
136	9.93488986041956e-07\\
137	9.58720398972191e-07\\
138	9.10961995419821e-07\\
139	8.9278826239423e-07\\
140	8.64196715661912e-07\\
141	8.40010748790737e-07\\
142	7.67012406483287e-07\\
143	7.40817527576935e-07\\
144	7.14401178887607e-07\\
145	6.91374762382663e-07\\
146	6.60948057770541e-07\\
147	6.42524320179846e-07\\
148	6.05050522514684e-07\\
149	5.84282609473051e-07\\
150	5.52897953684986e-07\\
151	5.15722193889432e-07\\
152	5.02297192957415e-07\\
153	4.97048170047204e-07\\
154	4.67395226471331e-07\\
155	4.58003456018318e-07\\
156	4.3457918985464e-07\\
157	4.17257844436837e-07\\
158	3.95760704408647e-07\\
159	3.87469138421331e-07\\
160	3.78458220029864e-07\\
161	3.70372300621757e-07\\
162	3.678520235106e-07\\
163	3.50082362173745e-07\\
164	3.40515686120067e-07\\
165	3.20951486905683e-07\\
166	3.05328926853698e-07\\
167	2.88620161196447e-07\\
168	2.78139573076759e-07\\
169	2.70718723978667e-07\\
170	2.61524838495892e-07\\
171	2.51934365102545e-07\\
172	2.50817495232691e-07\\
173	2.4571212113964e-07\\
174	2.26867444076557e-07\\
175	2.2205813839436e-07\\
176	2.13918777806436e-07\\
177	1.95612992358839e-07\\
178	1.87769452780772e-07\\
179	1.78620974202976e-07\\
180	1.69566034716707e-07\\
181	1.65688742977568e-07\\
182	1.64085510817964e-07\\
183	1.56769662006607e-07\\
184	1.50818686588491e-07\\
185	1.44032082767191e-07\\
186	1.42158390796685e-07\\
187	1.36500428084467e-07\\
188	1.31018438099678e-07\\
189	1.24034952020765e-07\\
190	1.19791420427383e-07\\
191	1.15615091133164e-07\\
192	1.1070055727196e-07\\
193	1.07872879763976e-07\\
194	1.05149051389442e-07\\
195	9.98487111201183e-08\\
196	9.79179411395028e-08\\
197	9.08048621095968e-08\\
198	8.87772089832042e-08\\
199	8.48748366991443e-08\\
200	8.19328752774286e-08\\
201	8.12056511547814e-08\\
202	8.08312628665413e-08\\
203	7.83480360245979e-08\\
204	7.52893032856099e-08\\
205	7.48322595975671e-08\\
206	7.34802206040113e-08\\
207	7.27824579862234e-08\\
208	6.71075009103122e-08\\
209	6.54070452715249e-08\\
210	6.24208804267449e-08\\
211	6.12539266608176e-08\\
212	6.112496614889e-08\\
213	5.85702001827922e-08\\
214	5.67371912015675e-08\\
215	5.47571541443672e-08\\
216	5.42227750548491e-08\\
217	5.17091779126089e-08\\
218	5.05224408699237e-08\\
219	4.93213895725901e-08\\
220	4.66445882833464e-08\\
221	4.63626261168661e-08\\
222	4.5896075933637e-08\\
223	4.56802948675963e-08\\
224	4.41044319891573e-08\\
225	4.11328543045964e-08\\
226	4.06950054747638e-08\\
227	3.71959335712709e-08\\
228	3.69417297003112e-08\\
229	3.66808407180334e-08\\
230	3.57302994134848e-08\\
231	3.45592253966229e-08\\
232	3.37295271275965e-08\\
233	3.26530835466142e-08\\
234	2.94860627767359e-08\\
235	2.91131854989663e-08\\
236	2.84911839975728e-08\\
237	2.81042118726182e-08\\
238	2.7907025955351e-08\\
239	2.6669194782553e-08\\
240	2.65770154680349e-08\\
241	2.6282907863117e-08\\
242	2.52030785792847e-08\\
243	2.4129460470081e-08\\
244	2.32258524289501e-08\\
245	2.28617443964005e-08\\
246	2.259166298289e-08\\
247	2.2437566847405e-08\\
248	2.12875594340831e-08\\
249	2.0635472876275e-08\\
250	1.89959131351329e-08\\
251	1.88509074085367e-08\\
252	1.86533791393392e-08\\
253	1.80714356471284e-08\\
254	1.7867877338182e-08\\
255	1.75536029691788e-08\\
256	1.70278888949544e-08\\
257	1.57794812860022e-08\\
258	1.5534055416378e-08\\
259	1.54332475655734e-08\\
260	1.50945600223088e-08\\
261	1.46604851199204e-08\\
262	1.41991648990659e-08\\
263	1.3273458513313e-08\\
264	1.26908930333067e-08\\
265	1.26612125745618e-08\\
266	1.23115333797232e-08\\
267	1.19414215069461e-08\\
268	1.16885258825765e-08\\
269	1.1518879567048e-08\\
270	1.133453154445e-08\\
271	1.13301595695368e-08\\
272	1.07787340833938e-08\\
273	1.06317932033846e-08\\
274	1.02200860161903e-08\\
275	1.01052191881889e-08\\
276	1.00165618489757e-08\\
277	9.63452987710929e-09\\
278	9.53961101852021e-09\\
279	8.87975340870222e-09\\
280	8.61103352886086e-09\\
281	8.49982400543906e-09\\
282	8.45396767734413e-09\\
283	8.0557925809559e-09\\
284	7.65228933928638e-09\\
285	7.42326251286923e-09\\
286	7.30529573620113e-09\\
287	7.09026239238552e-09\\
288	7.06115709344714e-09\\
289	6.71559107704281e-09\\
290	6.51804479688827e-09\\
291	6.43905330332042e-09\\
292	6.41299593468039e-09\\
293	6.162892751297e-09\\
294	6.09768900052467e-09\\
295	6.02298145652711e-09\\
296	5.9097158771351e-09\\
297	5.78726454215938e-09\\
298	5.35167618943006e-09\\
299	5.25829648537086e-09\\
300	5.13413528208122e-09\\
};
\end{axis}
\end{tikzpicture}%

%% file: fig/Pex2/param_err.tikz
\begin{tikzpicture}

\begin{axis}[%
width=\fwidth,
height=0.999\fheight,
at={(0\fwidth,0\fheight)},
scale only axis,
xmin=0.5,
xmax=1.1,
xlabel style={font=\color{white!15!black}},
xlabel={Parameters set},
ymode=log,
ymin=1e-05,
ymax=0.01,
yminorticks=true,
ylabel style={font=\color{white!15!black}},
ylabel={Relative error},
axis background/.style={fill=white},
xmajorgrids,
ymajorgrids,
yminorgrids,
legend style={legend cell align=left, align=left, draw=white!15!black}
]
\addplot [color=black, line width=2.0pt, only marks, draw=none, mark size=2.5pt, mark=o, mark options={solid, black}]
  table[row sep=crcr]{%
0.589048622548086	0.00128723606144595\\
0.719948316447661	0.000768654062005064\\
0.850848010347236	0.000694683832920954\\
0.98174770424681	0.000713239813809081\\
};
\addlegendentry{Test Parameters}

\addplot [color=black, line width=2.0pt, draw=none, only marks, mark size=1.8pt, mark=square, mark options={solid, black}]
  table[row sep=crcr]{%
0.523598775598299	0.00124901086611805\\
0.654498469497873	0.000861794091725826\\
0.785398163397448	0.000726077501233657\\
0.916297857297023	0.000693987103886093\\
1.0471975511966	0.000734535499491196\\
};
\addlegendentry{Training Parameters}

\addplot [color=red, line width=2.0pt, draw=none, mark size=2.5pt, mark=o, mark options={solid, red}, forget plot]
  table[row sep=crcr]{%
0.589048622548086	0.000410644878880693\\
0.719948316447661	0.000140646160393673\\
0.850848010347236	7.88018968226637e-05\\
0.98174770424681	8.11188155006359e-05\\
};
\addplot [color=red, line width=2.0pt, draw=none, mark size=1.8pt, mark=square, mark options={solid, red}, forget plot]
  table[row sep=crcr]{%
0.523598775598299	7.31370604531129e-05\\
0.654498469497873	9.15242745978392e-05\\
0.785398163397448	7.51169476943852e-05\\
0.916297857297023	6.39114836291818e-05\\
1.0471975511966	6.34705195138329e-05\\
};
\end{axis}
\end{tikzpicture}%